\documentclass[preprint,10pt]{elsarticle}

\usepackage{mathrsfs,amsmath,amssymb}
\usepackage{mathrsfs,amsmath,amssymb}
\usepackage{amscd} 
\usepackage{relsize,amsmath} 
\newcommand{\lsum}{\mathlarger{\sum}}

\newcommand{\lint}{\mathlarger{\int}}
\newcommand{\Lint}{\mathlarger{\mathlarger{\int}}}

\makeatletter
\newcommand{\biggg}[1]{{\hbox{$\left#1\vbox to 20.5pt{}\right.\n@space$}}}
\newcommand{\Biggg}[1]{{\hbox{$\left#1\vbox to 23.5pt{}\right.\n@space$}}}
\newcommand{\bigggg}[1]{{\hbox{$\left#1\vbox to 26.5pt{}\right.\n@space$}}}
\newcommand{\Bigggg}[1]{{\hbox{$\left#1\vbox to 29.5pt{}\right.\n@space$}}}
\newcommand{\biggggg}[1]{{\hbox{$\left#1\vbox to 32.5pt{}\right.\n@space$}}}
\newcommand{\Biggggg}[1]{{\hbox{$\left#1\vbox to 35.5pt{}\right.\n@space$}}}
\newcommand{\bigggggg}[1]{{\hbox{$\left#1\vbox to 38.5pt{}\right.\n@space$}}}
\newcommand{\Bigggggg}[1]{{\hbox{$\left#1\vbox to 41.5pt{}\right.\n@space$}}}
\makeatother
\newcommand{\bigggl}{\mathopen\biggg}

\newcommand{\bigggr}{\mathclose\biggg}
\newcommand{\Bigggl}{\mathopen\Biggg}

\newcommand{\Bigggr}{\mathclose\Biggg}

\makeatletter
\@addtoreset{equation}{section}
\makeatother

\begin{document}

\newtheorem{thm}{Theorem}
\newtheorem{lem}[thm]{Lemma}
\newdefinition{rmk}{Remark}
\newproof{pf}{Proof}
\newproof{pot}{Proof of Theorem \ref{thm2}}

\begin{frontmatter}



\title{Decay properties of solutions 
toward a multiwave pattern 
to the Cauchy problem 
for the scalar conservation law 
with degenerate flux and viscosity}

\author[rvt]{Natsumi Yoshida}
\ead{14v00067@gst.ritsumei.ac.jp}

\address{BKC Research Organization of Social Sciences, 
Ritsumeikan University, Kusatsu, Shiga 525-8577, Japan
/Osaka City University Advanced Mathematical Institute, 
Sumiyoshi, Osaka 558-8585, Japan.}
\address{}

\begin{abstract}
In this paper, we study the precise decay rate in time to solutions 
of the Cauchy problem for 
the one-dimensional conservation law 
with a nonlinearly degenerate viscosity 
$
\partial_tu +\partial_x \bigl(f(u) \bigr)
  = \mu \, 
    \partial_x \left( \, 
    \left| \, \partial_xu \, \right|^{p-1} \partial_xu \, 
    \right)
$
where the far field states are prescribed.
Especially, we deal with the case 
when the flux function is convex or concave 
but linearly degenerate on some interval. 
As the corresponding Riemann problem admits a Riemann solution as 
a multiwave pattern which consists 
of the rarefaction waves and the contact discontinuity, 
it has already been proved by Yoshida 
that the solution to the Cauchy problem 
tends toward the linear combination of 
the rarefaction waves and contact wave for $p$-Laplacian type viscosity 
as the time goes to infinity. 
We investigate that the decay rate in time 
of the corresponding solutions toward the multiwave pattern. 
Furthermore, we also investigate that the decay rate in time 
of the solution for the higher order derivative. 
The proof is given by 
$L^1$, $L^2$-energy and time-weighted $L^q$-energy methods 
under the use of the precise asymptotic properties of 
the interactions between the nonlinear waves. 
\end{abstract}

\end{frontmatter}
%




\pagestyle{myheadings}
\thispagestyle{plain}
\markboth{N. YOSHIDA}{NATSUMI YOSHIDA}

\section{Introduction and main theorems}
In this paper, 
we are concerned with the asymptotic behavior 
and the time-decay of 
solutions 
for the one-dimensional scalar conservation law 
with a nonlinearly degenerate viscosity 
($p$-Laplacian type viscosity with $p>1$) 
\begin{eqnarray}
 \left\{\begin{array}{ll}
  \partial_tu +\partial_x \bigl(f(u) \bigr)
  = \mu \, 
    \partial_x \left( \, 
    \left| \, \partial_xu \, \right|^{p-1} \partial_xu \, 
    \right)
  \qquad &(t>0, x\in \mathbb{R}), \\[5pt]
  u(0,x) = u_0(x) \qquad &( x \in \mathbb{R} ).
 \end{array}
 \right.\,
\end{eqnarray}
Here, $u=u(t,x)$ denotes the unknown function of $t>0$ and $x\in \mathbb{R}$, 
the so-called conserved quantity, 
$f=f(u)$ is the flux function depending only on $u$, 
$\mu$ is the viscosity coefficient, 
$u_0$ is the given initial data 
and $u$ is assumed to be 
asymptotically far field constant states $u_{\pm }$ 
as 
$
x\rightarrow \pm \infty,
$
that is, 
$$
\displaystyle{\lim_{x\to \pm \infty}} u(t,x) =u_{\pm}  \in \mathbb{R} 
\qquad \bigl( t \ge 0 \bigr).
$$
We suppose the given flux $f=f(u)$ is a $C^1$-function 
satisfying $f(0)=f'(0)=0$, 
$\mu$ is a positive constant 
and far field states $u_{\pm }$ satisfy $u_-<u_+$ 
without loss of generality. 

We are interested in the asymptotic behavior 
and its precise estimates in time of 
the global solution to the Cauchy problem (1.1). 
It can be expected that the large-time behavior 
is closely related to the weak solution (``Riemann solution'') of the corresponding 
Riemann problem 
(cf. \cite{lax}, \cite{smoller})
for the non-viscous hyperbolic part of (1.1): 
\begin{eqnarray}
 \left\{\begin{array} {ll}
 \partial _t u + \partial _x \bigl( f(u) \bigr)=0 
 \qquad &(t>0, x\in \mathbb{R}),\\[5pt]
u(0,x)=u_0 ^{\rm{R}} (x)\qquad &(x \in \mathbb{R}),
 \end{array}
  \right.\,
\end{eqnarray}
where $u_0 ^{\rm{R}}$ is the Riemann data defined by
$$
u_0 ^{\rm{R}} (x)=u_0 ^{\rm{R}} (x\: ;\: u_- ,u_+)
          := 
          \left\{\begin{array} {ll}
          u_-  & \; (x < 0),\\[5pt]
          u_+  & \; (x > 0).
          \end{array}\right.
$$
In fact, for $p=1$ in (1.1), the usual linear viscosity case: 
\begin{eqnarray}
 \left\{\begin{array}{ll}
  \partial_tu +\partial_x \bigl(f(u) \bigr)= \mu \, \partial_x^2 u
  \qquad &(t>0, x\in \mathbb{R}), \\[5pt]
  u(0,x) = u_0(x) \qquad &( x \in \mathbb{R} ),\\[5pt]
  \displaystyle{\lim_{x\to \pm \infty}} u(t,x) =u_{\pm}  
  \qquad &\bigl( t \ge 0 \bigr),   
 \end{array}
 \right.\,
\end{eqnarray}
when the smooth flux function $f$ is genuinely nonlinear 
on the whole space $\mathbb{R}$,
i.e., $f''(u)\ne 0\ (u\in \mathbb{R})$, 
Il'in-Ole{\u\i}nik \cite{ilin-oleinik} showed 
the following: if  $f''(u) > 0\ (u\in \mathbb{R})$, that is, the Riemann solution 
consists of a single rarefaction wave solution, 
the global solution in time of the Cauchy problem 
(1.3) tends toward the rarefaction wave;
if  $f''(u) < 0\ (u\in \mathbb{R})$, that is,
the Riemann solution consists of a single shock wave solution, 
the global solution of the
Cauchy problem (1.3) does the corresponding smooth traveling wave solution 
(``viscous shock wave'') of (1.3) with a spacial shift 
(cf. \cite{ilin-kalashnikov-oleinik}). 
Hattori-Nishihara \cite{hattori-nishihara} also proved that the asymptotic decay rate in time, 
of the solution toward the single rarefaction wave
(see also \cite{hashimoto-kawashima-ueda}, \cite{hashimoto-matsumura}, \cite{nakamura}). 
More generally, in the case of the flux functions 
which are not uniformly genuinely nonlinear, 
when the Riemann solution consists of a single shock wave 
satisfying Ole{\u\i}nik's shock condition, 
Matsumura-Nishihara \cite{matsu-nishi3} showed the asymptotic stability 
of the corresponding viscous shock wave. 
Moreover, Matsumura-Yoshida \cite{matsumura-yoshida} 
considered the circumstances 
where the Riemann solution generically forms 
a complex pattern of multiple nonlinear waves which consists of 
rarefaction waves 
and waves of contact discontinuity (refer to \cite{liep-rosh}), 
and investigated 
that the case where the flux function $f$ is smooth
and genuinely nonlinear (that is, $f$ is convex function or concave function) 
on the whole $\mathbb{R}$ except 
a finite interval $I := (a,b) \subset \mathbb{R}$, and 
linearly degenerate on $I$, that is, 
\begin{equation}
\left\{
\begin{array}{ll}
  f''(u) >0 & \; \bigl(u \in (-\infty ,a\, ]\cup [\, b,+\infty )\bigr),\\[5pt]
  f''(u) =0 & \; \bigl(u \in (a,b)\bigr).
\end{array}\right.
\end{equation}
Under the conditions (1.4), 
the corresponding Riemann solution does form 
multiwave pattern which consists of the contact discontinuity 
with the jump from $u=a$ to $u=b$ and the rarefaction waves, 
depending on the choice of $a$, $b$, $u_-$ and $u_+$.
Thanks to that the cases in which the interval $(a,b)$ is disjoint 
from the interval $(u_-,u_+)$ 
are similar as in the case 
the flux function $f$ is genuinely nonlinear 
on the whole space $\mathbb{R}$, 
and the case $u_-<a<u_+<b$ is 
the same as that for $a<u_-<b<u_+$, 
we may only consider the 
typical cases 
\begin{equation}
a<u_-<b<u_+ \quad \mbox{or} \quad u_-<a<b<u_+.   
\end{equation}
Under the conditions (1.4) and (1.5), 
they have shown the unique global solution in time 
to the Cauchy problem (1.3) tends uniformly in space toward 
the multiwave pattern of the combination of 
the viscous contact wave 
and the rarefaction waves 
as the time goes to infinity. 
It should be noted that the rarefaction wave which connects the 
far field states $u_-$ and $u_+$ 
$\bigl(u_\pm \in (-\infty,a\, ]$ or $u_\pm \in [\, b,\infty)\bigr)$
is explicitly given by
\begin{equation}
u=u^r \left( \frac{x}{t}\: ;\:  u_- , u_+ \right)
:= \left\{
\begin{array}{ll}
  u_-  & \; \bigl(\, x \leq \lambda(u_-)\,t \, \bigr),\\[7pt]
  \displaystyle{ (\lambda)^{-1}\left( \frac{x}{t}\right) } 
  & \; \bigl(\, \lambda(u_-)\,t \leq x \leq \lambda(u_+)\,t\,  \bigr),\\[7pt]
   u_+ & \; \bigl(\, x \geq \lambda(u_+)\,t \, \bigr),
\end{array}
\right.
\end{equation} 
where $\lambda(u):=f'(u)$,
and the viscous contact wave 
which connects 
$u_-$ and $u_+$ ($u_\pm \in [\,a,b\,]$) is given by an exact
solution of the linear 
convective heat equation 
\begin{equation}
  \partial _t u + \tilde{\lambda} \, \partial_x u = \mu \, \partial_x^2 u
  \qquad \biggl(
  \,\tilde{\lambda}:= \frac{f(b)-f(a)}{b-a},\ 
  t>0,x \in  \mathbb{R} 
  \biggr)
\end{equation}
 which has the form
$$
u=U\left(\frac{x-\tilde{\lambda} \, t}{\sqrt{t}}\: ;\: u_- ,u_+ \right)
$$
where $U\left(\frac{x}{\sqrt{t}}\,;\,u_- ,u_+ \right)$ is explicitly defined by
\begin{equation}
 U\left(\frac{x}{\sqrt{t}}\: ;\: u_- ,u_+ \right)
 :=u_- +\frac{u_+ - u_-}{\sqrt{\pi}}
   \lint ^{\frac{\mathlarger{x}}{\mathlarger{\sqrt{4\mu t}}}}_{-\infty} 
 \mathrm{e}^{-\xi^2}\, \mathrm{d}\xi 
\quad \, (t>0, x\in \mathbb{R}). 
\end{equation}
Yoshida \cite{yoshida} also obtained 
the almost optimal decay properties 
for the asymptotics toward the multiwave pattern. 
In fact, owing to \cite{yoshida}, 
the decay rate in time is 
$(1+t)^{-{\frac{1}{2}}\, \left(\frac{1}{2}-\frac{1}{p} \right)}$ 
in the $L^p$-norm $\bigl(2\leq p<+\infty \bigr)$ 
and 
$(1+t)^{-{\frac{1}{4}}\, +\epsilon}$ for any $\epsilon >0$ 
in the $L^\infty$-norm 
if the initial perturbation from the corresponding asymptotics 
satisfies $H^1$. 
Furthermore, 
if the perturbation satisfies $H^1\cap L^1$, 
the decay rate in time is 
$(1+t)^{-{\frac{1}{2}}\, \left(1-\frac{1}{p} \right)\, +\epsilon }$ 
for any $\epsilon >0$ 
in the $L^p$-norm $\bigl(1\leq p<+\infty \bigr)$ 
and
$(1+t)^{-{\frac{1}{2}}\, +\epsilon}$ for any $\epsilon >0$
in the $L^\infty$-norm. 

For $p>1$, 
there are few results for the asymptotic behavior 
for the problem (1.1)
(the related problems are studied in 
\cite{gur-mac}, \cite{nagai-mimura1}, \cite{nagai-mimura2} 
and so on). 
In the case where the flux function is genuinely nonlinear 
on the whole space $\mathbb{R}$, 
Matsumura-Nishihara \cite{matsu-nishi3} 
proved that if the far field states satisfy $u_- = u_+=:\tilde{u}$, 
then the solution tends toward the constant state $\tilde{u}$, 
and if the far field states $u_- < u_+$, 
then the solution tends toward a single rarefaction wave. 
Furthermore, Yoshida \cite{yoshida''} recently showed 
their precise decay estimates. 
In the case where the flux function is given as (1.4) 
and the far field states as (1.5), 
Yoshida \cite{yoshida'} also showed 
that the similar asymptotics as the one in \cite{matsumura-yoshida} 
which tends toward the multiwave pattern 
of the combination of 
the rarefaction waves 
which connect the 
far field states $u_-$ and $u_+$ 
$\bigl(u_\pm \in (-\infty,a\, ]$ or $u_\pm \in [\, b,\infty)\bigr)$, 
and the viscous contact wave which connects 
$u_-$ and $u_+$ ($u_\pm \in [\,a,b\,]$).  
In more detail, 
the viscous contact wave is 
said to be ``contact wave for $p$-Laplacian type viscosity'' 
and explicitly given by 
\begin{align}
U\left(\, \frac{x-\tilde{\lambda} \, t}
          {t^{\frac{1}{p+1}}}\: ;\: u_- ,u_+ \right) 
= u_- + \displaystyle{
        \lint^{\frac{\mathlarger{x-\tilde{\lambda} \, t}}
        { \mathlarger{t}^{{\scriptscriptstyle \frac{1}{p+1}}} }}_{-\infty}
             \Bigl( \left( \, 
                A-B \xi^2 \, 
                \right)\vee 0 \, \Bigr)^{\frac{1}{p-1}}
                \, \mathrm{d}\xi },
\end{align}
\begin{eqnarray*}
\left\{\begin{array}{ll}
\displaystyle{
A=A_{p,\mu,u_{\pm}}
:= \left( \, 
   \frac{(p-1)\, \left( \, u_{+} - u_{-} \, \right)}
   {8\, \mu\, p(p+1) 
   \biggl( \, \displaystyle{\int^{\frac{\pi}{2}}_{0} 
   \bigl( \sin \theta \bigr)^{\frac{p+1}{p-1}}\, \mathrm{d}\theta }
   \biggr)^2 }
   \, \right)^{\frac{p-1}{p+1}}
}, \\[25pt]
\displaystyle{ 
B=B_{p,\mu}:= \frac{p-1}{2\, \mu \, p(p+1)} 
}, \\[10pt] 
\displaystyle{ 
2A^{\frac{p+1}{2(p-1)}}B^{-\frac{1}{2}}
      \int^{\frac{\pi}{2}}_{0} 
      \bigl( \sin \theta \bigr)^{\frac{p+1}{p-1}}\, \mathrm{d}\theta 
      = u_{+} - u_{-}}.
\end{array}
\right.\,
\end{eqnarray*}
It should be noted that the wave (1.9) is 
constructed by the Barenblatt-Kompanceec-Zel'dovi{\v{c}} solution 
(see also \cite{carillo-toscani}, \cite{huang-pan-wang}, \cite{kamin}) 
\begin{equation}
v(t,x)
:= \frac{1}{(1+t)^{\frac{1}{p+1}}}\, 
   \Bigggl( \left( 
   A-B \left(\frac{x}{(1+t)^{\frac{1}{p+1}}} \right)^2 \, 
   \right)\vee 0 \Bigggr)^{\frac{1}{p-1}}
\end{equation}
of the following Cauchy problem of the porous medium equation 
\begin{eqnarray}
 \left\{\begin{array}{ll}
  \partial_t v
  = 
    \mu \, 
    \partial_x^2 
    \left( \, 
    \left| \, v \, \right|^{p-1} v \, 
    \right)
  \qquad &(t>-1, x\in \mathbb{R}), \\[5pt]
  v(-1,x) = (u_{+} - u_{-})\, \delta(x) 
  \qquad &( x \in \mathbb{R}\, ;\,  u_{-}<u_{+}),\\[5pt]
  \displaystyle{\lim_{x\to \pm \infty}} v(t,x) =0  
  \qquad &\bigl( t \ge -1 \bigr),  
 \end{array}
 \right.\,
\end{eqnarray}
where $\delta(x)$ is the Dirac $\delta$-distribution. 


\medskip
The aim of the present paper is to obtain 
the precise time-decay estimates for the asympotics of the 
previous study in \cite{yoshida'}. 

\medskip

\noindent
{\bf Stability Theorem } (Yoshida \cite{yoshida'}){\bf .}\quad{\it
Let the flux function $f$ satisfy {\rm(1.4)} and 
the far field states $u_{\pm }$ {\rm(1.5)}. 
Assume that the initial data satisfies 
$u_0-u_0 ^{\rm{R}} \in L^2$ and 
$\partial _x u_0 \in L^{p+1}$. 
Then the Cauchy problem {\rm(1.1)} with $p>1$ has a 
unique global 
solution in time $u=u(t,x)$ 
satisfying 
\begin{eqnarray*}
\left\{\begin{array}{ll}
u-u_0 ^{\rm{R}} \in C^0\bigl( \, [\, 0,\infty)\, ;L^2 \bigr)
                \cap L^{\infty}\bigl( \, \mathbb{R}^{+} \, ;L^{2} \bigr),\\[5pt]
\partial _x u \in L^{\infty} \bigl( \, \mathbb{R}^{+} \, ;L^{p+1} \bigr),\\[5pt]
\partial _t u 
\in L^{\infty} \bigl( \, \mathbb{R}^{+} \, ;L^{p+1} \bigr)
,\\[5pt]
\partial_x \left( \, \left| \, \partial_x u \, \right|^{p-1} \partial_xu \, \right)
\in L^{2}\bigl(\, {\mathbb{R}^{+}_{t}} \times {\mathbb{R}}_{x} \bigr),
\end{array} 
\right.\,
\end{eqnarray*}
and the asymptotic behavior 
$$
\lim _{t \to \infty}\sup_{x\in \mathbb{R}}
|\,u(t,x)-U_{multi}(\, t,x \: ;\: u_-,u_+)\,| = 0, 
$$
where $U_{multi}(t,x)=U_{multi}(\, t,x \: ;\: u_-,u_+)$ is defined as follows: 
in the case $a<u_-<b<u_+$,
$$
U_{multi}(t,x):=
U\left(\, \frac{x-\tilde{\lambda} \, t}{t^{\frac{1}{p+1}}}\: ;\: u_- , b \right)
 + u^r\left(\, \frac{x}{t}\: ;\: b , u_+ \right) - b
$$
and, in the case $u_-<a<b<u_+$,
$$
U_{multi}(t,x):=
u^r\left(\, \frac{x}{t}\: ;\: u_- ,a \right) -a
 + U\left(\, \frac{x-\tilde{\lambda} \, t}{t^{\frac{1}{p+1}}}\: ;\: a ,b \right)
 + u^r\left(\, \frac{x}{t}\: ;\: b , u_+ \right) - b.
$$
}

Now we are ready to state our main results. 

\medskip

\noindent
{\bf Theorem 1.1} (Main Theorem I){\bf .}\quad{\it
Under the same assumptions in Stability Theorem, 
the unique global solution in time $u$ 
of the Cauchy problem {\rm(1.1)} 
satisfying 
\begin{eqnarray*}
\left\{\begin{array}{ll}
u-u_0 ^{\rm{R}} \in C^0\bigl( \, [\, 0,\infty)\, ;L^2 \bigr)
                \cap L^{\infty}\bigl( \, \mathbb{R}^{+} \, ;L^{2} \bigr),\\[5pt]
\partial _x u \in L^{\infty} \bigl( \, \mathbb{R}^{+} \, ;L^{p+1} \bigr),\\[5pt]
\partial _t u 
\in L^{\infty} \bigl( \, \mathbb{R}^{+} \, ;L^{p+1} \bigr)
,\\[5pt]
\partial_x \left( \, \left| \, \partial_x u \, \right|^{p-1} \partial_xu \, \right)
\in L^{2}\bigl(\, {\mathbb{R}^{+}_{t}} \times {\mathbb{R}}_{x} \bigr),
\end{array} 
\right.\,
\end{eqnarray*}
satisfies the following time-decay estimates 
\begin{eqnarray*}
\left\{\begin{array} {ll}
\left|\left|\,
u(t) - 
U_{multi}(\, t,\cdot \: ;\: u_-,u_+)\,
\right|\right|_{L^q }
\leq C( \, p, \, q, \, u_0 \, ) \, 
     (1+t)^{-\frac{1}{3p+1}\left(1-\frac{2}{q}\right)},\\[7pt]
\left|\left|\,
u(t) - 
U_{multi}(\, t,\cdot \: ;\: u_-,u_+)
\,\right|\right|_{L^{\infty} }
\leq C( \, \epsilon, \, p, \, q, \, u_0, \, \partial _x u_0 \, ) \, 
     (1+t)^{-\frac{1}{3p+1}+\epsilon},
\end{array}
  \right.\,
\end{eqnarray*}
for 
$q \in [\, 2, \infty) $ and any $\epsilon>0$. 
}

\medskip

\noindent
{\bf Theorem 1.2} (Main Theorem I\hspace{-.1em}I){\bf .}\quad{\it
Under the same assumptions in Theorem 1.1, 
if the initial data further satisfies $u_0-u_0 ^{\rm{R}} \in L^1$, 
then it holds that the unique global solution in time $u$ 
of the Cauchy problem {\rm(1.1)} 
satisfies the following time-decay estimates 
\begin{eqnarray*}
\left\{\begin{array} {ll}
\left|\left| \,
u(t) - 
U_{multi}(\, t,\cdot \: ;\: u_-,u_+)
\, \right|\right|_{L^q }
\leq C( \, p, \, q, \, u_0 \, ) \, 
     (1+t)^{-\frac{1}{2p}\left(1-\frac{1}{q}\right)},\\[7pt]
\left|\left|\,
u(t) - 
U_{multi}(\, t,\cdot \: ;\: u_-,u_+)
\,\right|\right|_{L^{\infty} }
\leq C( \, \epsilon, \, p, \, q, \, u_0, \, \partial _x u_0 \, ) \, 
     (1+t)^{-\frac{1}{2p}+\epsilon},
\end{array}
  \right.\,
\end{eqnarray*}
for 
$q \in [\, 1, \infty) $ and any $\epsilon>0$. 
Furthermore, the solution satisfies 
the following time-decay estimates for the higher order derivative
\begin{align*}
\begin{aligned}
&\bigl|\bigl|\,
 \partial _x u(t) \,
 \bigr|\bigr|_{L^{p+1} }, \, \; \; 
 \left|\left|\,
 \partial _x u(t) - 
 \displaystyle{\partial _x U_{multi}(\, t,\cdot \: ;\: u_-,u_+)}
 \,\right|\right|_{L^{p+1} } \\[5pt]
& \leq \left\{\begin{array} {ll} 
      C( \, p, \, u_0, \, \partial _x u_0 \, ) \, 
      (1+t)^{-\frac{p}{(p+1)^2}} \, \: \: \; \; \; \quad \qquad 
      \left( \, 1 < p < 
      \displaystyle{\frac{7 + \sqrt{73} }{12} }
      \, \right),\\[15pt] 
      C( \, \epsilon, \, p, \, u_0, \, \partial _x u_0 \, ) \, 
      (1+t)^{-\frac{3}{2(p+1)(3p-2)} + \epsilon} \, \, \: \; \; \; 
      \left( \,  
      \displaystyle{\frac{7 + \sqrt{73} }{12} }
      \leq p \, \right),
      \end{array}
      \right.\,
\end{aligned}
\end{align*}
for any $\epsilon>0$. 
}

\medskip

\noindent
{\bf Theorem 1.3} (Main Theorem I\hspace{-.1em}I\hspace{-.1em}I ){\bf .}\quad{\it
Under the same assumptions in Theorem 1.2, 
if the initial data further satisfies 
$\partial _x u_0 \in L^{r+1} \, (r>p)$, 
then it holds that the unique global solution in time $u$ 
of the Cauchy problem {\rm(1.1)} 
satisfies the following time-decay estimates 
for the higher order derivative
\begin{align*}
\begin{aligned}
&\bigl|\bigl|\,
 \partial _x u(t) \,
 \bigr|\bigr|_{L^{r+1} }, \, \; \; 
 \left|\left|\,
 \partial _x u(t) - 
 \displaystyle{\partial _x U_{multi}(\, t,\cdot \: ;\: u_-,u_+)}
 \,\right|\right|_{L^{r+1} } \\[5pt]
& \leq \left\{\begin{array} {ll} 
       C( \, p, \, r, \, u_0, \, \partial _x u_0 \, ) \, 
       (1+t)^{-\frac{4p(r-p)+7p+3}{6p(p+1)(r+1)}} \\[15pt] 
       \, \, \, \: \; \; \; \quad \quad \quad \quad \qquad 
       \left( \, 
       1 < p < 
       \displaystyle{\frac{7 + \sqrt{73} }{12} }, \; 
       r > \displaystyle{\frac{-4p^2+7p+3}{2p}} > p 
       \, \right),\\[15pt] 
       C( \, \epsilon, \, p, \, r, \, u_0, \, \partial _x u_0 \, ) \, 
       (1+t)^{-\frac{r}{(p+1)(r+1)}+ \epsilon} \\[15pt] 
       \, \, \: \; \; \; \quad \quad \quad \quad \quad \qquad 
       \left( \, 
       1 < p < 
       \displaystyle{\frac{7 + \sqrt{73} }{12} }, \; 
       p < r \le \displaystyle{\frac{-4p^2+7p+3}{2p}}
       \, \right),\\[15pt] 
       C( \, \epsilon, \, p, \, r, \, u_0, \, \partial _x u_0 \, ) \, 
       (1+t)^{-\frac{p+2r}{2p(3p-2)(r+1)}+ \epsilon} \\[15pt] 
       \, \, \, \: \: \; \; \quad \quad \qquad \qquad \qquad \qquad \qquad 
       \left( \, 
       \displaystyle{\frac{7 + \sqrt{73} }{12}} \leq p
       \, \right),
       \end{array}
       \right.\,
\end{aligned}
\end{align*}
for any $\epsilon>0$. 
}

\medskip

This paper is organized as follows. 
In Section 2, we shall prepare the basic properties of 
the rarefaction wave and the 
contact wave for $p$-Laplacian type viscosity.  
In Section 3, 
we reformulate the problem in terms of the deviation from 
the asymptotic state 
(similarly in \cite{matsumura-yoshida}, \cite{yoshida}, \cite{yoshida'}). 
Following the arguments in \cite{matsu-nishi2}, 
we also prepare some uniform boundedness and energy estimates 
of the deviation 
as the solution to the reformulated problem. 
We further introduce the precise properties of the interactions 
between the nonlinear waves. 
In order to obtain the time-decay estimates (Theorem 1.1 and Theorem 1.2), 
in Section 4 and Section 5, 
we establish the uniform energy estimates in time 
by using a very technical time-weighted energy method. 
Finally, in Section 6, 
we prove the time-decay $L^{r+1}$-estimate 
for the higher order derivative, Theorem 1.3. 

\smallskip

{\bf Some Notation.}\quad 
We denote by $C$ generic positive constants unless 
they need to be distinguished. 
In particular, use $C(\alpha, \beta, \cdots )$ 
or $C_{\alpha, \beta, \cdots }$ 
when we emphasize the dependency on $\alpha, \beta, \cdots $. 
Use $\mathbb{R}^{+}$ as 
$
\mathbb{R}^{+}:=(0,\infty),
$
and the symbol ``$\vee $'' as
$$
a\vee b:= \max \{a,b\}. 
$$
We also use the Friedrichs mollifier $\rho_\delta \ast $, 
where, 
$\rho_\delta(x):=\frac{1}{\delta}\rho \left( \frac{x}{\delta}\right)$
with 
\begin{align*}
\begin{aligned}
&\rho \in C^{\infty}_0(\mathbb{R}),\quad 
\rho (x)\geq 0\:  (x \in \mathbb{R}), \\
&\mathrm{supp} \{\rho \} \subset 
\left\{x \in \mathbb{R}\: \left|\:  |\, x \, |\le 1 \right. \right\},\quad  
\int ^{\infty}_{-\infty} \rho (x)\, \mathrm{d}x=1, 
\end{aligned}
\end{align*}
and $\rho_\delta \ast f$ denote the convolution. 
For function spaces, 
$L^p = L^p(\mathbb{R})$ and $H^k = H^k(\mathbb{R})$ 
denote the usual Lebesgue space and 
$k$-th order Sobolev space on the whole space $\mathbb{R}$ 
with norms $||\cdot||_{L^p}$ and $||\cdot||_{H^k}$, respectively. 
We also define 
the bounded $C^{m}$-class $\mathscr{B}^{m}$ as follows 
$$
f\in \mathscr{B}^{m}(\Omega)
\, \Longleftrightarrow  \, 
f\in C^m(\Omega), 
\; 
\sup _{\Omega}\, \sum _{k=0}^{m} \, \bigl| \, D^kf\, \bigr|<\infty 
$$
for $m< \infty$ and 
$$
f\in \mathscr{B}^{\infty }(\Omega)
\, \Longleftrightarrow  \, 
\forall n\in \mathbb{N},\, f\in C^n(\Omega), 
\; 
\sup _{\Omega}\, \sum _{k=0}^{n} \, \bigl| \, D^kf\, \bigr|<\infty 
$$
where $\Omega \subset \mathbb{R}^d$ and 
$D^k$ denote the all of $k$-th order derivatives. 

\section{Preliminaries} 
In this section, 
we shall arrange the several lemmas concerned with 
the basic properties of 
the rarefaction wave 
for accomplishing the proof of our main theorems. 
Since the rarefaction wave $u^r$ is not smooth enough, 
we need some smooth approximated one 
as in the previous results in 
\cite{hashimoto-matsumura}, \cite{liu-matsumura-nishihara}, 
\cite{matsu-nishi1}, \cite{matsumura-yoshida}, 
\cite{yoshida}, \cite{yoshida'}, \cite{yoshida''}. 
We start with the well-known arguments on $u^r$ 
and the method of constructing its smooth approximation. 
We first consider the rarefaction wave solution $w^r$ 
to the Riemann problem 
for the non-viscous Burgers equation 
\begin{equation}
\label{riemann-burgers}
  \left\{\begin{array}{l}
  \partial _t w + 
  \displaystyle{ \partial _x \left( \, \frac{1}{2} \, w^2 \right) } = 0 
  \, \, \; \; \qquad \quad \qquad ( t > 0,\,x\in \mathbb{R}),\\[7pt]
  w(0,x) = w_0 ^{\rm{R}} ( \, x\: ;\: w_- ,w_+):= \left\{\begin{array}{ll}
                                                  w_+ & \, \: \; \quad (x>0),\\[5pt]
                                                  w_- & \, \: \; \quad (x<0),
                                                  \end{array}
                                                  \right.
  \end{array}
  \right.
\end{equation}
where $w_\pm \in \mathbb{R}$ are 
the prescribed far field states. 
If the far field states satisfy $w_-<w_+$, 
then the Riemann problem (2.1) has a unique global weak solution 
$w=w^r\left( \, \frac{x}{t}\: ;\: w_-,w_+\right)$ 
given explicitly by 
\begin{equation}
\label{rarefaction-burgers}
w^r \left( \, \frac{x}{t}\: ;\: w_-,w_+\right) := 
  \left\{\begin{array}{ll}
  w_{-} & \bigl(\, x \leq w_{-}t \, \bigr),\\[5pt]
  \displaystyle{ \frac{x}{t} } & \bigl(\, w_{-} t \leq x \leq w_{+} t \, \bigr),\\[5pt]
  w_+ & \bigl(\, x\geq w_{+} t \, \bigr).
  \end{array}\right.
\end{equation} 
Next, under the condition 
$f''(u)>0\ (u\in \mathbb{R})$ and $u_-<u_+$, 
the rarefaction wave solution 
$u=u^r\left( \, \frac{x}{t}\: ;\: u_-,u_+\right)$ 
of the Riemann problem (1.2) 
for hyperbolic conservation law 
is exactly given by 
\begin{equation}
u^r\left( \, \frac{x}{t} \: ; \:  u_-,u_+\right) 
= (\lambda)^{-1}\biggl( w^r\left( \, \frac{x}{t} \: ; \:  \lambda_-,\lambda_+\right)\biggr)
\end{equation}
which is nothing but (1.6), 
where $\lambda_\pm := \lambda(u_\pm) = f'(u_\pm)$. 
We also consider 
the Cauchy problem for the following 
non-viscous Burgers equation
\begin{eqnarray}
\label{smoothappm}
\left\{\begin{array}{l}
 \partial _t w 
 + \displaystyle{ \partial _x \left( \, \frac{1}{2} \, w^2 \right) } =0 
 \, \, \; \; \quad \qquad \qquad \qquad \qquad \qquad  
 (\ t>0,\,x\in \mathbb{R}),\\[7pt]
 w(0,x) 
 = w_0(x) 
 := \displaystyle{ \frac{w_-+w_+}{2} + \frac{w_+-w_-}{2}\tanh x } 
 \qquad \quad \; \:  (x\in \mathbb{R}). 
\end{array}
\right.
\end{eqnarray}   
By using the method of characteristics, 
we easily see that the Cauchy problem (2.4) 
has a unique classical solution 
$$
w=w(\, t,x\: ;\: w_-,w_+)\in \mathscr{B}^{\infty }( \, [\, 0,\infty )\times \mathbb{R})
$$
satisfying the following formula
\begin{eqnarray}
 \left\{\begin{array} {l}
 w(t,x)=w_0\bigl(x_0(t,x)\bigr)=
 \displaystyle{ \frac{\lambda_-+\lambda_+}{2} } 
+ \displaystyle{ \frac{\lambda_+-\lambda_-}{2}\tanh \bigl( x_0(t,x) \bigr)} ,\\[7pt]
 x=x_0(t,x)+w_0\bigl(x_0(t,x)\bigr)\,t.
 \end{array}
  \right.\,
\end{eqnarray}
We define a smooth approximation of $w^r(\, \frac{x}{t}\: ;\: w_-,w_+)$ 
by the solution $w$. 
We also note the assumption of the flux function $f$ to be 
$\lambda'(u)\left( =\frac{\mathrm{d}^2f}{\mathrm{d}u^2}(u)\right)>0$. 

Now we summarize the results for the smooth approximation $w(\, t,x\: ;\: w_-,w_+)$ 
in the next lemma. 
Since the proof is given by the direct calculation as in \cite{matsu-nishi1}, 
we omit it. 

\medskip

\noindent
{\bf Lemma 2.1.}\quad{\it
Assume that the far field states satisfy $w_-<w_+$. 
Then the classical solution $w(t,x)=w(\, t,x\: ;\: w_-,w_+)$
given by {\rm(2.4)} 
satisfies the following properties: 

\noindent
{\rm (1)}\ \ $w_- < w(t,x) < w_+$ and\ \ $\partial_xw(t,x) > 0$  
\quad  $(t>0, x\in \mathbb{R})$.

\smallskip

\noindent
{\rm (2)}\ For any $1\leq q \leq \infty$, there exists a positive 
constant $C_q$ such that
             \begin{eqnarray*}
                 \begin{array}{l}
                    \parallel \partial_x w(t)\parallel_{L^q} \leq 
                    C_q (1+t)^{-1+\frac{1}{q}} 
                    \; \quad \bigl(t\ge 0 \bigr),\\[5pt]
                    \parallel \partial_x^2 w(t) \parallel_{L^q} \leq 
                    C_q (1+t)^{-1} 
                    \, \, \quad \quad \bigl(t\ge 0 \bigr).
                    \end{array}       
              \end{eqnarray*}
              
\smallskip

\noindent
{\rm (3)}\ $\: \displaystyle{\lim_{t\to \infty} 
\sup_{x\in \mathbb{R}}
\left| \,w(t,x)- w^r \left( \frac{x}{t} \right) \, \right| = 0}.$
}

\bigskip

\noindent
We define the approximation for 
the rarefaction wave $u^r\left( \, \frac{x}{t}\: ;\: u_-,u_+\right)$ by 
\begin{equation}
U^r(\, t,x\: ; \: u_-,u_+) := (\lambda)^{-1} \bigl( w(\, t,x\: ;\: \lambda_-,\lambda_+)\bigr).
\end{equation}

Then we have the next lemma 
as in the previous works 
(cf. \cite{hashimoto-matsumura}, \cite{liu-matsumura-nishihara}, 
\cite{matsu-nishi1}, \cite{matsumura-yoshida}, 
\cite{yoshida}, \cite{yoshida'}, \cite{yoshida''}). 

\medskip

\noindent
{\bf Lemma 2.2.}\quad{\it
Assume that the far field states satisfy $u_-<u_+$, 
and the flux fanction $f\in C^3(\mathbb{R})$, $f''(u)>0 \: (u\in [\,u_-,u_+\,])$. 
Then we have the following properties:

\noindent
{\rm (1)}\ $U^r(t,x)$ defined by {\rm (2.6)} is 
the unique $C^2$-global solution in space-time 
of the Cauchy problem
$$
\left\{
\begin{array}{l} 
\partial _t U^r +\partial _x \bigl( f(U^r ) \bigr) = 0 
\, \, \, \, \; \; \; \quad \quad \qquad \qquad \qquad \qquad  
(t>0, x\in \mathbb{R}),\\[7pt]
U^r(0,x) 
= \displaystyle{ (\lambda)^{-1} \left( \, \frac{\lambda_- + \lambda_+}{2} 
+ \frac{\lambda_+ - \lambda_-}{2} \tanh x \, \right) } 
\; \: \: \quad \quad( x\in \mathbb{R}),\\[7pt]
\displaystyle{\lim_{x\to \pm \infty}} U^r(t,x) =u_{\pm} 
\, \, \: \: \; \; \quad \quad \qquad \qquad \qquad \qquad \qquad \qquad 
\bigl(t\ge 0 \bigr).
\end{array}
\right.\,     
$$
{\rm (2)}\ \ $u_- < U^r(t,x) < u_+$ and\ \ $\partial_xU^r(t,x) > 0$  
\quad  $(t>0, x\in \mathbb{R})$.

\smallskip

\noindent
{\rm (3)}\ For any $1\leq q \leq \infty$, there exists a positive 
constant $C_q$ such that
             \begin{eqnarray*}
                 \begin{array}{l}
                    \parallel \partial_x 
                    U^r(t) 
                    \parallel_{L^q} \leq 
                    C_q(1+t)^{-1+\frac{1}{q}} 
                    \quad \bigl(t\ge 0 \bigr),\\[5pt]
                    \parallel \partial_x^2 U^r(t) \parallel_{L^q} \leq 
                    C_q(1+t)^{-1}
                    \, \, \quad \quad \bigl(t\ge 0 \bigr).
                    \end{array}       
              \end{eqnarray*}
              
\smallskip

\noindent
{\rm (4)}\ $\: \displaystyle{\lim_{t\to \infty} 
\sup_{x\in \mathbb{R}}
\left| \,U^r(t,x)- u^r \left( \frac{x}{t} \right) \, \right| = 0}.$

\smallskip

\noindent
{\rm (5)}\ For any $\epsilon \in (0,1)$, there exists a positive 
constant $C_\epsilon$ such that
$$
\left|\, 
U^r(t,x)-u_+ \, \right|
\leq C_\epsilon(1+t)^{-1+\epsilon}\mathrm{e}^{-\epsilon |x-\lambda_+t|}
\quad \bigl(t\ge 0, x \ge \lambda_+t \bigr).
$$

\smallskip

\noindent
{\rm (6)}\ For any $\epsilon \in (0,1)$, there exists a positive 
constant $C_\epsilon$ such that
$$
\left|\, 
U^r(t,x)-u_- \, \right|
\leq C_\epsilon(1+t)^{-1+\epsilon}\mathrm{e}^{-\epsilon |x-\lambda_-t|}
\quad \bigl(t\ge 0, x \le \lambda_-t\bigr).
$$

\smallskip

\noindent
{\rm (7)}\ For any $\epsilon \in (0,1)$, there exists a positive 
constant $C_\epsilon$ such that
$$
\left| \,U^r(t,x) - u^r\left( \frac{x}{t}\right) \, \right| 
\leq C_\epsilon(1+t)^{-1+\epsilon} 
\qquad \bigl(t \ge 1,  \lambda_-t \le x \le \lambda_+t\bigr).
$$

\smallskip

\noindent
{\rm (8)}\ For any $(\epsilon ,q)\in (0,1)\times [\, 1,\infty \, ]$, 
there exists a positive 
constant $C_{\epsilon,q}$ such that
$$
\left|\left|
 \,U^r(t,\cdot \: ) - u^r\left( \frac{\cdot }{t}\right) \, 
\right|\right|_{L^q}  
\leq C_{\epsilon,q}(1+t)^{-1+\frac{1}{q}+\epsilon} 
\qquad \bigl(t \ge 0 \bigr).
$$
}

\noindent
Because the proofs of (1) to (4) are given in \cite{matsu-nishi1}, 
(5) to (7) are in \cite{matsumura-yoshida} 
and (8) is in \cite{yoshida}, 
we omit the proofs here. 
\bigskip

We also prepare the next lemma for the properties of 
the contact wave for $p$-Laplacian type viscosity 
$U\Bigl(\, \frac{x}{t^{\frac{1}{p+1}}}\,;\,u_- ,u_+ \Bigr)$ defined by (1.11). 
In the following, we abbreviate 
``contact wave for $p$-Laplacian type viscosity'' 
to ``viscous contact wave''.  
we rewrite the viscous contact wave as 
\begin{align}
\begin{aligned}
U(t,x)&=U\left(\, \frac{x}{t^{\frac{1}{p+1}}}\: ;\: u_- ,u_+ \right) \\
      &=u_{+}-\Lint^{\infty}_{x}
             \frac{1}{t^{\frac{1}{p+1}}}\, 
             \bigggl( \left( 
             A-B \left(\frac{y}{t^{\frac{1}{p+1}}} \right)^2 \, 
             \right)\vee 0 \bigggr)^{\frac{1}{p-1}}
             \, \mathrm{d}y, 
\end{aligned}
\end{align}
where 
\begin{eqnarray*}
\left\{\begin{array}{ll}
\displaystyle{
A=A_{p,\mu,u_{\pm}}
:= \left( \, 
   \frac{(p-1)\, \left( \, u_{+} - u_{-} \, \right)}
   {8\, \mu\, p(p+1) 
   \biggl( \, \displaystyle{\int^{\frac{\pi}{2}}_{0} 
   \bigl( \sin \theta \bigr)^{\frac{p+1}{p-1}}\, \mathrm{d}\theta }
   \biggr)^2 }
   \, \right)^{\frac{p-1}{p+1}}
}, \\[25pt]
\displaystyle{ 
B=B_{p,\mu}:= \frac{p-1}{2\, \mu \, p(p+1)} 
}, \\[10pt] 
\displaystyle{ 
2A^{\frac{p+1}{2(p-1)}}B^{-\frac{1}{2}}
      \int^{\frac{\pi}{2}}_{0} 
      \bigl( \sin \theta \bigr)^{\frac{p+1}{p-1}}\, \mathrm{d}\theta 
      = u_{+} - u_{-}}
.
\end{array}
\right.\,
\end{eqnarray*}

\noindent
Then, the following properties hold. 

\medskip

\noindent
{\bf Lemma 2.3.}\quad{\it
For any 
$p>1$ and $u_\pm \in \mathbb{R}$, we have the following: 

\noindent
{\rm (i)}\ $U$ defined by {\rm (1.11)} satisfies 
$$
U \in 
\mathscr{B}^1\bigl( \, (0,\infty \bigr)\times \mathbb{R}\, ) 
\mathlarger{ \mathlarger{ \mathlarger{ 
\setminus } } } \, 
C^2
\left( \, 
\biggl\{ \, 
(t,x) \in \mathbb{R}^{+} \times \mathbb{R} \, \biggr|\biggl. \, 
x = \pm \sqrt{\frac{A}{B}}\, t^{\frac{1}{p+1}} \, 
\biggr\} \, \right), 
$$  
and is a 
self-similar type strong solution of the Cauchy problem 
   $$
           \left\{
              \begin{array}{l} 
              \partial _t U 
              - \mu \, \partial_x \left( \, 
              \left| \, \partial_x U \, \right|^{p-1} \partial_x U \, \right) 
              = 0 
              \quad \qquad \qquad  (t>0, x\in \mathbb{R}),\\[13pt]
              U(0,x) = u_0 ^{\rm{R}} (\, x\: ;\: u_- ,u_+)
               = 
               \left\{\begin{array} {ll}
               u_-  & \, \; \; \quad \qquad (x < 0),\\[5pt]
               u_+  & \, \; \; \quad \qquad (x > 0),
               \end{array}\right.\\[13pt]
              \displaystyle{\lim_{x\to \pm \infty}} U(t,x) =u_{\pm} 
              \, \: \: \: \; \quad \quad \qquad \qquad \qquad \qquad \bigl(\, t\ge 0\bigr).
              \end{array}
            \right.\,   
    $$          
{\rm (ii)}\ For $t>0$ and $x\in \mathbb{R}$, 
$$
\left\{
   \begin{array}{l} 
   U(t,x)=u_{-},
   \, \: \; \quad \qquad \qquad \qquad \qquad 
   \left( \, x \leq -\sqrt{\frac{A}{B}}\, t^{\frac{1}{p+1}} \,  \right),
   \\[5pt]
   u_- < U(t,x) < u_+,\; \partial_xU(t,x) > 0 \quad 
   \left( \, 
   -\sqrt{\frac{A}{B}}\, t^{\frac{1}{p+1}} < x < \sqrt{\frac{A}{B}}\, t^{\frac{1}{p+1}}
   \, \right),
   \\[5pt]
   U(t,x)=u_{+},
   \, \: \; \quad \qquad \qquad \qquad \qquad 
   \left( \,  x \geq \sqrt{\frac{A}{B}}\, t^{\frac{1}{p+1}} \,  \right). 
   \end{array}
\right.\, 
$$

\smallskip

\noindent
{\rm (iii)}\ It holds that for any $1 \leq  q < \infty $, 
$$
\|\, \partial _x 
U (t)\, \| _{L^q} 
= C_1(\, A,B\: ;\: p ,q\, ) \, t^{-\frac{q-1}{(p+1)q}}\qquad (t > 0) 
$$
where 
$$
C_1(\, A,B\: ;\: p ,q\, ):=
\left( 
2A^{\frac{p+2q-1}{2(p-1)}}B^{-\frac{1}{2}}
\int^{\frac{\pi}{2}}_{0} 
      \bigl( \sin \theta \bigr)^{\frac{q}{p-1}}
\, \mathrm{d}\theta 
\right)^{\frac{1}{q}}.
$$
If $q = \infty$, we have 
$$
\|\, \partial _x U (t) \, \| _{L^\infty} 
= \left(2A \right)^{\frac{1}{p-1}} \, t^{- \frac{1}{p+1}}\qquad (t > 0).
$$

\smallskip

\noindent
{\rm (iv)}\ It holds that for any $1 \leq  q < \frac{p-1}{p-2}$ with $p>2$, 
or any $1 \leq  q < \infty$ with $1 < p \leq 2$, 
$$
\|\, \partial _x^2 U (t) \, \| _{L^q} 
= C_2(\, A,B\: ;\: p ,q\, ) \, t^{-\frac{2q-1}{(p+1)q}}\qquad (t > 0) 
$$
where 
\begin{align*}
\begin{aligned}
& C_2(\, A,B\: ;\: p ,q\, ) \\
& \quad \, 
:=
\bigggl( 
2\left(
\frac{2A^{-\frac{ p-2 }{p-1}} B}{p-1}
\right)^{q}
\left( \frac{B}{A} \right)^{- \frac{q+1}{2}}
\int^{\frac{\pi}{2}}_{0} 
      \bigl( \sin \theta \bigr)^{- \frac{2 ( p-2 )q}{p-1} + 1}
      \bigl( \cos \theta \bigr)^{q}
\, \mathrm{d}\theta
\bigggr)^{\frac{1}{q}}. 
\end{aligned}
\end{align*}
If $1 < p \leq 2$, for $q = \infty$, we have 
$$
\|\, \partial _x^2 U (t) \, \| _{L^\infty} 
= {\frac{2A^{\frac{\left| p-2 \right|}{p-1}} B}{p-1}} 
     \left( \frac{B}{A} \right)^{-\frac{1}{2}}\, 
     t^{- \frac{2}{p+1}}\qquad (t > 0).
$$

\smallskip

\noindent
{\rm (v)}\ It holds that 
$$
\left|\left| \, 
\partial_x \left( \, 
\left| \, \partial_x U \, \right|^{p-1} \partial_x U \, 
\right) (t) \, 
\right| \right| _{L^2} 
= C_3(\, A,B\: ;\: p \, ) \, t^{-\frac{2p+1}{2(p+1)}}\qquad (t > 0) 
$$
where 
$$
C_3(\, A,B\: ;\: p \, ):=
\left( 
2\left(
\frac{2B^{p}}{p-1}
\right)^{2}
\left( \frac{B}{A} \right)^{-\frac{3p-7}{2(p-1)}}
\int^{\frac{\pi}{2}}_{0} 
      \bigl( \sin \theta \bigr)^{\frac{p+3}{p-1}}
      \bigl( \cos \theta \bigr)^{2}
\, \mathrm{d}\theta
\right)^{\frac{1}{2}}.
$$
\smallskip

\noindent
{\rm (vi)}\; $\displaystyle{\lim_{t\to \infty} \sup_{x\in \mathbb{R}}\, 
\bigl|\,U(1+t,x)- U(t,x)\, \bigr| = 0}.$
             
\smallskip

\noindent
{\rm (vii)}\ For any $1\leq q \leq \infty$ with $p>1$, 
there exists a positive 
constant $C_{p,q}$ such that
$$
\left|\left|
 \,U(1+t,\cdot \: ) - U(t,\cdot \: ) \, 
\right|\right|_{L^q}  
\leq C_{p,q}\, t^{-1+\frac{1}{(p+1)q}} 
\qquad \bigl(t > 0 \bigr).
$$

}


\section{Reformulation of the problem}
In this section, 
we reduce our Cauchy problem (1.1) to 
a simpler case 
and reformulate the problem 
in terms of the deviation from the asymptotic state 
(the same as in \cite{matsumura-yoshida}, \cite{yoshida}, \cite{yoshida'}). 
At first, 
without loss of generality, 
we shall consider the case 
where $a<0$, $b=0$ and the flux function $f(u)$ satisfies 
\begin{equation}
\left\{
\begin{array}{ll}
  f''(u) >0 & \; (u \in (-\infty ,a]\cup [0,+\infty )),\\[5pt]
  f(u) =0 & \; (u \in (a,0)),
\end{array}\right.
\end{equation}
under changing the variables and constant as
$x-\tilde{\lambda} \,t \mapsto x$,  $u-b \mapsto u$,
$f(u+b)-f'(b)\,u-f(a) \mapsto f(u)$ and $a-b \mapsto a$
in this order. 
For the far field states 
$u_\pm \in \mathbb{R}$, 
we only deal with the typical case 
$a<u_-<0<u_+$ for simplicity, since the case 
$u_-<a<0<u_+$ can be treated 
technically 
in the same way of the proof as 
$a<u_-<0<u_+$. 
Indeed, in the case $u_-<a<0<u_+$,
as we shall see in Section 4  and Section 5, 
there appears the extra nonlinear interaction terms
between two rarefaction waves 
$u^r(\, \frac{x}{t}\: ;\: u_-,a)$ and $u^r(\, \frac{x}{t}\: ;\: 0,u_+)$
with $\lambda(a)=\lambda(0)=0$ in the remainder term of 
the viscous conservation law for the asymptotics 
$U_{multi}$\,(see the right-hand side of (3.4)). 
These terms can be handled in much easier way 
by Lemma 2.2 than that for other essential nonlinear interaction terms 
between the rarefaction and the viscous contact waves. 
Furthermore, we should point out that 
the problem under the assumptions for the flux function (3.1) 
and the far field states $a<u_-<0<u_+$
is essentially the same as that for $a=-\infty$, because 
obtaining the a priori and the uniform energy estimates for the former one can be
given in almost the same way as the latter one. 
Therefore, it is quite natural for us 
to treat only 
a simple case 
\begin{equation}
\left\{
\begin{array}{ll}
  f''(u) >0 & \; (u \in [\, 0,\infty )),\\[5pt]
  f(u) =0 & \; (u \in (-\infty,0)),
\end{array}\right.
\end{equation}
and assume $u_-<0<u_+$.
The corresponding stability theorem 
and our main theorems are the following. 

Under the situation, we reformulate the problem 
in terms of the deviation from the asymptotic state.  
We first should note by Lemma 2.2 and Lemma 2.3, 
the asymptotic state $u^r\left( \, \frac{x}{t}\: ;\: u_-,u_+\right)$ 
can be replaced by a following approximated multiwave pattern 
\begin{equation*}
\tilde{U}(t,x) := U(1+t,x) + U^r(t,x),
\end{equation*}
where
$$
U(1+t,x)=U\left( \frac{x}{(1+t)^{\frac{1}{p+1}}}\: ;\: u_- ,0 \right),
\quad U^r(t,x) =U^r(\, t,x\: ; \: 0, u_+).
$$
In fact, 
from Lemma 2.1 (especially (8)) and Lemma 2.3 (especially (vii)), 
it follows that 
for any $\epsilon>0$
\begin{align*}
\begin{aligned}
\left|\left| \, 
\tilde{U}(t,\cdot \: ) - U_{multi}(t,\cdot \: ) \, \right|\right|_{L^q}
&\le
\left|\left| \, U(1+t,\cdot \: ) - U(t,\cdot \: ) \, \right|\right|_{L^q}
+ \left|\left| \, U^r(t,\cdot \: ) 
- u^r\left( \frac{\cdot}{t}\right) \, \right|\right|_{L^q}\\
&\le C_{\epsilon,q}(1+t)^{-\left( 1-\frac{1}{q}\right)+\epsilon}
\qquad (t\ge 0\, ;\, 1\le q\le \infty).
\end{aligned}
\end{align*}
Then it is noted that $\tilde{U}$ is 
monotonically increasing and approximately satisfies 
the equation of (1.1) as
\begin{equation}
 \partial _t\tilde{U} +\partial_x \bigl(f(\tilde{U})\bigr)
 - \mu \, 
 \partial_x 
\left( \, 
\bigl| \, \partial_x \tilde{U} \, \bigr|^{p-1} \partial_x \tilde{U} \, 
\right)
 = -F_{p}(U,U^r),
\end{equation}
where the remainder term $F_{p}(U,U^r)$ is explicitly given by
\begin{align}
\begin{aligned}
F_{p}(U,U^r)
:= 
&\widetilde{F_{p} }(U,U^r)\\
&+ \mu \, \partial_x 
\left( \, 
\bigl| \, \partial_x U + \partial_x U^{r} \, \bigr|^{p-1} 
\bigl( \, \partial_x U + \partial_x U^{r} \, \bigr) - 
\left| \, \partial_x U \, \right|^{p-1} \partial_x U  \, 
\right)\\
:= 
&-\bigl( \, f'(U+U^r)
- f'(U^r) \, \bigr) \, \partial_x U^r 
- f'(U+U^r)\, \partial_x U \\
&+ \mu \, \partial_x 
\left( \, 
\bigl| \, \partial_x U + \partial_x U^{r} \, \bigr|^{p-1} 
\bigl( \, \partial_x U + \partial_x U^{r} \, \bigr) - 
\left| \, \partial_x U \, \right|^{p-1} \partial_x U  \, 
\right)
\end{aligned}
\end{align}
which consists of the interaction terms of the viscous contact wave $U$ and 
the approximation of the rarefaction wave $U^r$, 
and the approximation error of $U^r$ as solution to 
the conservation law for the $p$-Laplacian type viscosity. 
Here we should 
note that 
$U$ is monotonically nondecreasing 
and $U^r$ is monotonically increasing, 
that is, 
$\partial_x\tilde{U}(t,x) > 0\ \bigl(t\ge 0, x\in \mathbb{R} \bigr)$
which is frequently used hereinafter. 
Now  putting 
\begin{equation}
u(t,x) = \tilde{U}(t,x) + \phi(t,x)
\end{equation}
and using (3.5),
we can reformulate the problem (1.1) in terms of 
the deviation $\phi $ from 
$\tilde{U}$ as 
\begin{eqnarray}
 \left\{\begin{array}{ll}
  \partial _t\phi + \partial_x \left( f(\tilde{U}+\phi) - f(\tilde{U}) \right) \\[5pt]
  \quad - \mu \, \partial_x 
    \left( \, 
    \bigl| \, \partial_x \tilde{U} + \partial_x \phi \, \bigr|^{p-1} 
    \bigl( \, \partial_x \tilde{U} + \partial_x \phi \, \bigr) - 
    \bigl| \, \partial_x \tilde{U} \, \bigr|^{p-1} \partial_x \tilde{U}  \, 
    \right)\\[2pt]
    \quad \qquad \qquad \qquad \qquad \qquad \qquad = F_{p}(U,U^r)
  \, \, \: \: \quad  (t>0, x\in \mathbb{R}), \\[5pt]
  \phi(0,x) = \phi_0(x) 
  := u_0(x)-\tilde{U}(0,x) 
  \qquad \qquad \quad \; \: \: \; \; \: \; \; \,\, \, \, (x\in \mathbb{R}), \\[5pt]
  \displaystyle{\lim _{x\rightarrow \pm \infty}\, \phi (t,x) = 0 }
  \qquad \qquad \qquad \qquad \qquad \qquad \qquad \quad \, \bigl( t \geq 0 \bigr).
 \end{array}
 \right.\,
\end{eqnarray}
Then we look for the global solution in time 
$$
\phi \in C^0\bigl( \, [\, 0,\infty)\, ;L^2 \bigr)
         \cap L^{\infty}\bigl( \, \mathbb{R}^{+} \, ;L^{2} \bigr) 
$$ with 
$$
\partial _x \phi \in L^{\infty} \bigl( \, \mathbb{R}^{+} \, ;L^{p+1} \bigr)
                     \cap L^{p+1}
                     \bigl(\, {\mathbb{R}^{+}_{t}} 
                     \times {\mathbb{R}}_{x} \bigr). 
$$
Here we note the fact $\phi_0 \in L^2$ 
and $\partial_x \phi_0 \in L^{p+1}$ by the assumptions on $u_0$ 
and the fact 
$$
\partial_x \tilde{U}(0,\cdot \, )
=\partial_x U(0,\cdot \, )+\partial_x U^r(0,\cdot \, )\in L^{p+1}.
$$
In the following, we always assume that 
the flux function $f\in C^1(\mathbb{R})\cap C^3 (\, [\, 0,\infty))$ satisfies (3.2), 
and  
the far field states satisfy $u_-<0<u_+$. 
Then the corresponding our main theorems for $\phi$ 
we should prove are as follows. 

\medskip

\noindent
{\bf Theorem 3.1.}\quad{\it
Assume that the flux function 
$f\in C^1(\mathbb{R})\cap C^3 (\, [\, 0,\infty))$ satisfies {\rm(3.2)},  
the far field states $u_-<0<u_+$, 
and the initial data 
$\phi_0 \in L^2$ and 
$\partial _x u_0 \in L^{p+1}$. 
Then, the unique global solution in time $\phi$ 
of the Cauchy problem {\rm(3.6)} 
satisfying 
\begin{eqnarray*}
\left\{\begin{array}{ll}
\phi \in C^0\bigl( \, [\, 0,\infty)\, ;L^2 \bigr)
            \cap L^{\infty}\bigl( \, \mathbb{R}^{+} \, ;L^{2} \bigr),\\[5pt]
\partial _x \phi \in L^{\infty} \bigl( \, \mathbb{R}^{+} \, ;L^{p+1} \bigr)
                     \cap L^{p+1}
                          \bigl(\, {\mathbb{R}^{+}_{t}} 
                          \times {\mathbb{R}}_{x} \bigr),\\[5pt]
\partial _x \bigl( \, \tilde{U} + \phi \, \bigr)
              \in L^{\infty} \bigl( \, \mathbb{R}^{+} \, ;L^{p+1} \bigr)
                  \cap L^{p+2}\bigl(\, {\mathbb{R}^{+}_{t}} 
                       \times 
                       \left\{x \in \mathbb{R}\, |\, u>0 \right\} \bigr),\\[5pt]
\partial _t \bigl( \, \tilde{U} + \phi \, \bigr) 
\in  L^{\infty} \bigl( \, \mathbb{R}^{+} \, ;L^{p+1} \bigr)
,\\[5pt]
\partial_x \left( \, \bigl| \, 
\partial_x\bigl( \, \tilde{U} + \phi \, \bigr)\, 
\bigr|^{p-1} \partial_x\bigl( \, \tilde{U} + \phi \, \bigr) \, \right)
\in L^{2}\bigl(\, {\mathbb{R}^{+}_{t}} \times {\mathbb{R}}_{x} \bigr),
\end{array} 
\right.\,
\end{eqnarray*}
and 
$$
\displaystyle{\lim_{t\to \infty}}\sup_{x\in \mathbb{R}}\, 
\bigl| \, \phi(t,x)\,\bigr| = 0
$$
satisfies the following time-decay estimates 
\begin{eqnarray*}
\left\{\begin{array} {ll}
\left|\left|\,
\phi (t) 
\, \right|\right|_{L^q }
\leq C( \, p, \, q, \, \phi_0 \, ) \, 
     (1+t)^{-\frac{1}{3p+1}\left(1-\frac{2}{q}\right)},\\[7pt]
\left|\left|\,
\phi (t) 
\, \right|\right|_{L^{\infty} }
\leq C( \, \epsilon, \, p, \, q, \, \phi_0, \, \partial _x u_0 \, ) \, 
     (1+t)^{-\frac{1}{3p+1}+\epsilon},
\end{array}
  \right.\,
\end{eqnarray*}
for 
$q \in [\, 2, \infty) $ and any $\epsilon>0$. 
}

\medskip

\noindent
{\bf Theorem 3.2.}\quad{\it
Under the same assumptions in Theorem 3.1, 
if the initial data further satisfies $\phi_0 \in L^1$, 
then it holds that the unique global solution in time $\phi$ 
of the Cauchy problem {\rm(3.6)} 
satisfies the following time-decay estimates 
\begin{eqnarray*}
\left\{\begin{array} {ll}
\left|\left| \,
\phi (t) 
\, \right|\right|_{L^q }
\leq C( \, p, \, q, \, \phi_0 \, ) \, 
     (1+t)^{-\frac{1}{2p}\left(1-\frac{1}{q}\right)},\\[7pt]
\left|\left|\,
\phi (t) 
\,\right|\right|_{L^{\infty} }
\leq C( \, \epsilon, \, p, \, q, \, \phi_0, \, \partial _x u_0 \, ) \, 
     (1+t)^{-\frac{1}{2p}+\epsilon},
\end{array}
  \right.\,
\end{eqnarray*}
for 
$q \in [\, 1, \infty) $ and any $\epsilon>0$. 
Furthermore, the solution satisfies 
the following time-decay estimates for the higher order derivative
\begin{align*}
\begin{aligned}
&\bigl|\bigl|\,
 \partial _x u(t) \,
 \bigr|\bigr|_{L^{p+1} }, \, \; \; 
 \left|\left|\,
 \partial _x \phi (t) 
 \,\right|\right|_{L^{p+1} } \\[5pt]
& \leq \left\{\begin{array} {ll} 
      C( \, \epsilon, \, p, \, \phi_0, \, \partial _x u_0 \, ) \, 
      (1+t)^{-\frac{p}{(p+1)^2}} \, \quad 
      \left( \, 1 < p \le 
      \displaystyle{\frac{7}{12} + 
      \sqrt{\frac{73}{144} - \frac{(p+1)^2(3p-2)}{3}\, \epsilon } }
      \, \right),\\[15pt] 
      C( \, \epsilon, \, p, \, \phi_0, \, \partial _x u_0 \, ) \, 
      (1+t)^{-\frac{3}{2(p+1)(3p-2)} + \epsilon} 
      \left( \,  
      \displaystyle{\frac{7}{12} + 
      \sqrt{\frac{73}{144} - \frac{(p+1)^2(3p-2)}{3}\, \epsilon } }
      < p \, \right),
      \end{array}
      \right.\,
\end{aligned}
\end{align*}
for any $0<\epsilon \ll 1$. 
}

\medskip

\noindent
{\bf Theorem 3.3.}\quad{\it
Under the same assumptions in Theorem 3.2, 
if the initial data further satisfies 
$\partial _x u_0 \in L^{r+1} \, (r>p)$, 
then it holds that the unique global solution in time $\phi$ 
of the Cauchy problem {\rm(3.6)} 
satisfies the following time-decay estimates 
for the higher order derivative
\begin{align*}
\begin{aligned}
&\bigl|\bigl|\,
 \partial _x u(t) \,
 \bigr|\bigr|_{L^{r+1} }, \, \; \; 
 \left|\left|\,
 \partial _x \phi(t) 
 \,\right|\right|_{L^{r+1} } \\[5pt]
& \leq \left\{\begin{array} {ll} 
       C( \, \epsilon, \, p, \, r, \, \phi_0, \, \partial _x u_0 \, ) \, 
       (1+t)^{-\frac{4p(r-p)+7p+3}{6p(p+1)(r+1)}} \\[10pt] 
       \left( \, 
       1 < p \le 
       \displaystyle{\frac{7}{12} + 
       \sqrt{\frac{73}{144} - \frac{p(p+1)(3p-2)(r+1)}{2(r-p+1)}\, \epsilon } }, \; 
       r > \displaystyle{\frac{-4p^2+7p+3}{2p}} > p 
       \, \right),\\[15pt] 
       C( \, \epsilon, \, p, \, r, \, \phi_0, \, \partial _x u_0 \, ) \, 
       (1+t)^{-\frac{p+2r}{2p(3p-2)(r+1)}+ \epsilon} \\[15pt] 
       \left( \, 
       1 < p \le 
       \displaystyle{\frac{7}{12} + 
       \sqrt{\frac{73}{144} - \frac{p(p+1)(3p-2)(r+1)}{2(r-p+1)}\, \epsilon } }, \; 
       p < r \le \displaystyle{\frac{-4p^2+7p+3}{2p}}
       \, \right),\\[15pt] 
       C( \, \epsilon, \, p, \, r, \, \phi_0, \, \partial _x u_0 \, ) \, 
       (1+t)^{-\frac{6p(r-p)+7p+2r+3}{2(3p+1)(3p-2)(r+1)}+ \epsilon} \\[10pt] 
       \, \, \, \, \: \: \; \; \; \; \quad \quad \quad \quad \qquad \qquad \qquad 
       \left( \, 
       \displaystyle{\frac{7}{12} + 
       \sqrt{\frac{73}{144} - \frac{p(p+1)(3p-2)(r+1)}{2(r-p+1)}\, \epsilon } } < p
       \, \right),
       \end{array}
       \right.\,
\end{aligned}
\end{align*}
for any $0<\epsilon \ll 1$. 
}

\medskip
In order to accomplish the proofs of 
Theorem 3.1, Theorem 3.2 and Theorem 3.3, 
we will need some estimates about boundedness 
of the perturbation $\phi$ and $u$.  
We shall arrange some lemmas for them. 

By using the maximum principle 
(cf. \cite{ilin-kalashnikov-oleinik}, \cite{ilin-oleinik}), 
we first have the following uniform boundedness 
of the perturbation $\phi$ (and also $u$), that is, 

\medskip

\noindent
{\bf Lemma 3.1} (uniform boundedness){\bf .}\quad {\it
It holds that 
\begin{align}
\begin{aligned}
&\sup_{t\in [\, 0,\infty ), x\in \mathbb{R}}|\,\phi(t,x)\,| \\
&\le 
\sup_{x\in \mathbb{R}}|\,u_0 (x)\,| 
+ \sup_{t\in [\, 0,\infty ), x\in \mathbb{R}}|\,U (t,x)\,| 
+ \sup_{t\in [\, 0,\infty ), x\in \mathbb{R}}|\,U^r (t,x)\,| \, \\ 
&= 
\Vert \,\phi_0 \, \Vert_{L^{\infty}} 
+ 2 \, \bigl(\, | \, u_- \, | + | \, u_+ \, | \, \bigr), 
\end{aligned}
\end{align}
\begin{align}
\begin{aligned}
&\sup_{t\in [\, 0,\infty), x\in \mathbb{R}}|\, u(t,x)\,| \\
&\le \| \,\phi_0\, \|_{L^{\infty}} 
     + 2 \, \bigl(\, | \, u_- \, | + | \, u_+ \, | \, \bigr)
     + | \, u_- \, | \vee | \, u_+ \, | 
     =: \widetilde{C}. 
     \quad \quad \; \; \: \: \, \, 
\end{aligned}
\end{align}
}

\medskip

\noindent
Secondly, 
we also have the uniform estimates of $\phi$ as follows
(for the proof of it, 
see \cite{yoshida'}). 

\medskip

\noindent
{\bf Lemma 3.2} (uniform estimates){\bf .}\quad {\it
The unique global solution in time $\phi$ of the Cauchy problem {\rm (3.3)} 
satisfies the following uniform energy inequalities. 

\smallskip

\noindent
{\rm (1)}\ 
There exists a positive constant 
$
C_{p}(\phi_0)=C_{p}\bigl(\, \| \, \phi _0 \, \| _{L^2} \, \bigr)
$
such that 
\begin{align*}
\begin{aligned}
&\| \, \phi (t)\, \| _{L^2}^2
+\int _0^{\infty} G(t) \, \mathrm{d}t \\
&+\int _0^{\infty} \int _{-\infty}^{\infty} 
\bigl( \partial _x \phi \bigr)^2 
\left( 
\bigl| \partial _x \phi \bigr|^{p-1} 
+ \bigl| \partial _x U \bigr|^{p-1} + \bigl|  \partial _x U^r  \bigr|^{p-1}  
\right) \, \mathrm{d}x \mathrm{d}t 
\leq C_{p}(\phi_0)
 \quad \bigl(t \ge 0\bigr),
\end{aligned}
\end{align*}
where $G=G(t)$ is exactly given by 
\begin{align*}
\begin{aligned}
G(t):=\left( \, \int_{ \tilde{U} \geq 0} 
     \phi^2 \, \partial _x \tilde{U}\, \mathrm{d}x\right) (t)   
     &+\left( \, \int_{\tilde{U} + \phi \geq 0, \tilde{U} < 0} 
     \bigl( \, \tilde{U}+\phi \, \bigr)^2 \partial _x \tilde{U}
     \, \mathrm{d}x\right) (t)\\
    &+\left( \, \int_{\tilde{U} + \phi < 0, \tilde{U} \geq 0} 
     \bigl( \, \tilde{U}+|\, \phi \, |\,  \bigr)^2 \partial _x \tilde{U}
     \, \mathrm{d}x\right) (t).
\end{aligned}
\end{align*}
\smallskip

\noindent
{\rm (2)}\ 
There exists a positive constant 
$
C_{p}(\phi_0, \partial _x u_0)
=C_{p}\bigl(\, \| \, \phi _0 \, \| _{L^2}, 
 \| \, \partial _x u_0 \, \| _{L^{p+1}} \, \bigr)
$
such that 
\begin{align*}
\begin{aligned}
\| \, \partial _x u (t) 
\, \| _{L^{p+1}}^{p+1} 
&+\int _0^{\infty} \int _{-\infty}^{\infty} 
 \bigl| \, \partial _x u \, \bigr|^{2(p-1)} 
 \left( \, \partial _x^2 u \, \right)^2 
 \, \mathrm{d}x \mathrm{d}t \\
&\qquad \quad \leq 
 C_{p}(\phi_0, \partial _x u_0)
  \quad \bigl(t \ge 0\bigr).
\end{aligned}
\end{align*}
\smallskip

\noindent
{\rm (3)}\ 
There exists a positive constant 
$
C_{p}(\phi_0, \partial _x u_0)
=C_{p}\bigl(\, \| \, \phi _0 \, \| _{L^2}, 
 \| \, \partial _x u_0 \, \| _{L^{p+1}} \, \bigr)
$
such that 
\begin{equation*}
\int _0^{\infty}
 \| \, 
 \partial _x u (t) 
 \, \| _{L^{p+2}\left( \left\{x \in \mathbb{R}\, |\, u>0 \right\} \right)}^{p+2}
 \, \mathrm{d}t 
\leq C_{p}(\phi_0, \partial _x u_0)
 \quad \bigl( t \ge 0 \bigr).
\end{equation*}
}

\medskip

We also prepare the precise properties 
for the nonlinear interaction terms 
of the viscous contact wave $U$ 
and the approximation of rarefaction wave $U^r$. 
Namely, due to Lemma 2.2 and Lemma 2.3, 
we can easily see the fact that 
for any $t\ge 0$ there uniquely exists $x=X(t)\in \mathbb{R}$ such that 
\begin{equation}
\tilde{U}\bigl(t,X(t)\bigr) 
= U\bigl(t,X(t)\bigr)+U^r\bigl(t,X(t)\bigr)= 0
\quad \bigl( t \ge 0 \bigr).
\end{equation}
that is,
\begin{align*}
\begin{aligned}
U^r\bigl(t,X(t)\bigr)
&= -U\bigl(t,X(t)\bigr) \\ 
&= \Lint^{\infty}_{X(t)}
  \frac{1}{(1+t)^{\frac{1}{p+1}}}\, 
  \bigggl( \left( 
  A-B \left(\frac{y}{(1+t)^{\frac{1}{p+1}}} \right)^2 
  \right)\vee 0 \bigggr)^{\frac{1}{p-1}}
  \, \mathrm{d}y \quad \bigl( t \ge 0 \bigr). 
\end{aligned}
\end{align*}
More precisely, we have the following lemma. 

\medskip

\noindent
{\bf Lemma 3.3.}\quad {\it
The function 
$$
X:[\, 0,\infty)\ni t \longmapsto X(t)\in \mathbb{R}
$$
defined by {\rm(3.9)} 
has following asymptotic properties. 
\smallskip

\noindent
{\rm (i)}\ 
There exists a positive time $T_0$ such that 
for some $\delta \in \left(\, 0, \sqrt{\frac{A}{B}} \, \right)$, 
\begin{equation*}
\left(\, \sqrt{\frac{A}{B}} - \delta \, \right)\, (1+t)^{\frac{1}{p+1}} 
< X(t) < \sqrt{\frac{A}{B}}\, (1+t)^{\frac{1}{p+1}}
\qquad \bigl(\, t \ge T_0 \, \bigr).
\end{equation*}
\smallskip

\noindent
{\rm (ii)}\
For any $\epsilon \in (0,1)$, 
there exists a positive constant $C_{p, \epsilon}$ such that 
\begin{equation*}
\left|\, 
(\lambda)^{-1} \left(\frac{X(t)}{1+t} \right) -
\lint^{\infty}_{\frac{\mathlarger{X(t)}}{ \mathlarger{(1+t)^{\frac{1}{p+1}} }}}
                \Bigl ( \left( \, 
                A-B \, \xi^2 \, 
                \right)\vee 0 \Bigr)^{\frac{1}{p-1}}
                \, \mathrm{d}\xi 
\, \right| \le C_{p, \epsilon}(1+t)^{-1 + \epsilon} 
\qquad \bigl(\, t \ge T_0 \, \bigr).
\end{equation*}
\smallskip

\noindent
{\rm (iii)}\ 
There exists a positive constant $C_{p}$ such that 
$$
\left| \, 
\sqrt{\frac{A}{B}} - \frac{X(t)}{ (1+t)^{\frac{1}{p+1}} } 
\, \right|
\leq 
     C_{p} (1+t)^{-{\frac{p-1}{p+1}}}
\qquad \bigl(\, t \ge T_0 \, \bigr).
$$
}

\bigskip 

\noindent
\section{Time-decay estimates with $2 \le q \le \infty$ }
In this section, we show the time-decay estimates with $2 \le q \le \infty$ 
(not assuming $L^1$-integrability to the initial perturbation), 
that is, Theorem 3.1. 
To do that, we shall obtain the time-weighted $L^q$-energy estimates 
to $\phi$ with $2 \le q< \infty$ (cf. \cite{yoshida}). 

\medskip

\noindent
{\bf Proposition 4.1.}\quad {\it
Suppose the same assumptions in Theorem 3.1. 
For any $q \in [ \, 2, \infty)$, 
there exist positive constants $\alpha$ and $C_{\alpha,p,q}$, 
such that the unique global solution in time $\phi$ of the Cauchy problem {\rm(3.6)} 
satisfies the following $L^q$-energy estimate 
 \begin{align}
 \begin{aligned}
 &(1+t)^\alpha  \| \,\phi(t) \,\|_{L^q}^q 
  + \int ^t_0 (1+\tau )^\alpha  G_{q}(\tau) 
    \, \mathrm{d}\tau \\
 & \quad 
  + \int ^t_0 (1+\tau )^\alpha  \int _{-\infty}^{\infty} 
    | \, \phi \, |^{q-2} 
    \bigl( \, \partial _x \phi \, \bigr)^2 \\
 & \qquad \qquad \qquad \qquad \quad \; \, 
    \times 
    \left( \, 
    \bigl| \partial _x \phi \bigr|^{p-1} 
    + \bigl| \partial _x U \bigr|^{p-1} 
    + \bigl| \partial _x U^r \bigr|^{p-1}  
    \, \right) \, \mathrm{d}x \mathrm{d}\tau  \\
 & \quad 
  + \int ^t_0 (1+\tau )^\alpha  \int _{-\infty}^{\infty} 
    | \, \phi \, |^{q-2} \, 
    \left| \, 
    \bigl| \partial _x \phi + \partial _x U + \partial _x U^r \bigr|^{p-1} 
    - \bigl| \partial _x U + \partial _x U^r \bigr|^{p-1}  
    \, \right| \\
 & \qquad \qquad \qquad \qquad \quad \; \, 
   \times 
    \left| \, 
    \bigl( \, \partial _x \phi + \partial _x U + \partial _x U^r \, \bigr)^2 
    - \bigl( \, \partial _x U + \partial _x U^r \, \bigr)^2 \, \right| 
    \, \mathrm{d}x \mathrm{d}\tau  \\
 &\leq C_{\alpha,p,q} \| \,\phi_0 \,\|_{L^q}^q 
       + C\left( \, \alpha, \, p, \, q, \, \phi_0 \, \right) \, 
       (1+t)^{\alpha - \frac{q-2}{3p+1}}
       \quad \bigl(\, t \ge T_0 \, \bigr),
 \end{aligned}
 \end{align}
 where $G_q=G_q(t)$ is explicitly given by 
\begin{align*}
\begin{aligned}
G_q(t)&:=\left( \, \int_{\tilde{U}+\phi \geq 0, \tilde{U} \geq 0} 
        | \, \phi \, |^q \partial _x \tilde{U}\, \mathrm{d}x \, \right) (t) \\
      &\, \quad \;   +\left( \, \int_{\tilde{U}+\phi < 0, \tilde{U} \geq 0}
         | \, \phi \, |^{q-1}\, \tilde{U} \, 
         \partial _x \tilde{U} \, \mathrm{d}x \, \right) (t)
         + \left( \, \int_{\tilde{U}+\phi < 0,\tilde{U} \geq 0}
           {\tilde{U}}^q\, \partial _x \tilde{U} \, \mathrm{d}x \, \right) (t)\\
      &\, \quad \;  +\left( \int_{\tilde{U} + \phi \geq 0, \tilde{U} < 0} 
        \left( \, | \, \phi \, |^{q-1} 
        \bigl(\, q \, \tilde{U} + (q-1)\, | \, \phi \, | \, \bigr) + | \, \tilde{U} \, |^q\right)
        \partial _x \tilde{U}\, \mathrm{d}x \, \right) (t). \\
\end{aligned}
\end{align*}
}

\medskip

The proof of Proposition 4.1 is provided by the
following two lemmas.

\medskip

\noindent
{\bf Lemma 4.1.}\quad {\it
For any $2\leq q < \infty$, 
there exist positive constants $\alpha$ and $C_q$ such that 
\begin{align}
\begin{aligned}
&(1+t)^\alpha  \| \,\phi(t) \,\|_{L^q}^q 
  + q \, (q-1)
    \int ^t_0 (1+\tau )^\alpha  G_{q}(\tau) 
    \, \mathrm{d}\tau \\
 & 
  + C_{q} 
    \int ^t_0 (1+\tau )^\alpha 
    \int _{-\infty}^{\infty} 
    | \, \phi \, |^{q-2} 
    \bigl( \, \partial _x \phi \, \bigr)^2 \\
 & \qquad \qquad \qquad \qquad \quad \quad \quad \; \, 
   \times 
    \left( \, 
    \bigl| \partial _x \phi \bigr|^{p-1} 
    + \bigl| \partial _x U \bigr|^{p-1} 
    + \bigl| \partial _x U^r \bigr|^{p-1} 
    \, \right) \, \mathrm{d}x \mathrm{d}\tau  \\
 & 
  + C_{q} 
    \int ^t_0 (1+\tau )^\alpha  \int _{-\infty}^{\infty} 
    | \, \phi \, |^{q-2} \, 
    \left| \, 
    \bigl| \partial _x \phi + \partial _x \tilde{U} \bigr|^{p-1} 
    - \bigl| \partial _x \tilde{U} \bigr|^{p-1}  
    \, \right| \\
 & \qquad \qquad \qquad \qquad \quad \quad \quad \; \, 
   \times 
    \left| \, 
    \bigl( \, \partial _x \phi + \partial _x \tilde{U} \, \bigr)^2 
    - \bigl( \, \partial _x \tilde{U} \, \bigr)^2 \, \right| 
    \, \mathrm{d}x \mathrm{d}\tau  \\
 &\leq \| \,\phi_0 \,\|_{L^q}^q 
       + \alpha 
       \int ^t_0 (1+\tau )^{\alpha -1} 
       \| \,\phi(\tau ) \,\|_{L^p}^p 
       \, \mathrm{d}\tau \\
 &\, \: \quad \quad \quad \quad \, + 
    q \int ^t_0 (1+\tau )^\alpha 
    \| \,\phi(\tau ) \,\|_{L^\infty}^{p-1} 
    \left|\left|\, \widetilde{F_{p} }(U,U^r)(\tau ) \, \right| \right|_{L^1}
    \, \mathrm{d}\tau \\
 &\, \: \quad \quad \quad \quad \, + 
    \mu \, q \int ^t_0 (1+\tau )^\alpha 
    \| \,\phi(\tau ) \,\|_{L^\infty}^{p-2} \\
 &\, \: \quad \quad \quad \quad \quad \, \times    
    \int _{-\infty}^{\infty} 
    \bigl| \, \partial _x \phi \, \bigr|
    \Bigl( \, 
    \bigl( \, \partial_x U + \partial_x U^{r} \, \bigr)^{p} 
    - \left( \, \partial_x U \, \right)^{p} \, 
     \Bigr) 
     \, \mathrm{d}x\mathrm{d}\tau \quad \bigl( t \ge 0 \bigr). 
 \end{aligned}
 \end{align}
}

\bigskip

\noindent
{\bf Lemma 4.2.}\quad {\it
Assume $p>1$ and $2\leq q < \infty$. 
We have the following interpolation inequalities.  

 \noindent
 {\rm (1)}\ \ For any $2\leq r < \infty$, there exists a positive 
 constant $C_{p,q,r}$  such that 
 \begin{align*}
 \begin{aligned}
 &\Vert \, \phi (t) \, \Vert _{L^r }
  \leq 
  C_{p,q,r} \left( \, \int _{-\infty}^{\infty} | \, \phi \, |^2 \, \mathrm{d}x \, \right)
  ^{\frac{pr+p+q-1}{(3p+q-1)r}} \\
 & \qquad \quad \quad \quad \quad \; \: \, 
         \times 
         \left( \, \int _{-\infty}^{\infty} | \, \phi \, |^{q-2} 
         \bigl| \, \partial _x \phi \, \bigr|^{p+1} \, \mathrm{d}x \, \right)
         ^{\frac{r-2}{(3p+q-1)r}} \quad \bigl( t \ge 0 \bigr). 
 \end{aligned}
 \end{align*}

 \noindent
 {\rm (2)}\ \ 
 There exists a positive 
 constant $C_{p,q}$ 
 such that 
 \begin{align*}
 \begin{aligned}
 &\Vert \, \phi (t)  \, \Vert _{L^\infty }
  \leq 
  C_{p,q} \left( \, \int _{-\infty}^{\infty} | \, \phi \, |^2 \, \mathrm{d}x \, \right)
  ^{\frac{p}{3p+q-1}} \\
 & \qquad \quad \quad \quad \quad \; \, 
         \times 
         \left( \, \int _{-\infty}^{\infty} | \, \phi \, |^{q-2} 
         \bigl| \, \partial _x \phi \, \bigr|^{p+1} \, \mathrm{d}x \, \right)
         ^{\frac{1}{3p+q-1}} \quad \bigl( t \ge 0 \bigr). 
 \end{aligned}
 \end{align*}
}

\medskip

In what follows, we first prove Lemma 4.1 
and finally give the proof of Proposition 4.1 
(the proof of Lemma 4.2 is given in \cite{yoshida''}, 
so we omit here). 

\medskip

{\bf Proof of Lemma 4.1}.
Multiplying the equation in (3.6) by 
$\left|\phi \right|^{q-2} \phi$ with $2\leq q < \infty$, 
we obtain the divergence form 
 \begin{align}
 \begin{aligned}
 &\partial_t\left(\frac{1}{q} \left|\, \phi \, \right|^q \right) \\
 &+\partial _x 
   \biggl( \, 
   \left|\, \phi \, \right|^{q-2} \phi \, 
   \left( f(\tilde{U}+\phi)-f(\tilde{U}) \right) 
   \biggr) \\ 
 &+\partial _x \left( 
   -(q-1)\int _0^{\phi} 
   \left( f(\tilde{U}+\eta)-f(\tilde{U}) \right)\left|
   \, \eta \, \right|^{q-2}
   \, \mathrm{d}\eta 
   \, \right) \\
 &+\partial _x \biggl( \, 
   -\mu \, \left|\, \phi \, \right|^{q-2} \phi \, 
   \biggr.\\
 & \quad \times \biggl. 
    \Bigl( \, 
    \bigl|\, \partial _x \tilde{U} + \partial _x \phi \, \bigr|^{p-1} 
    \bigl( \, \partial _x \tilde{U} + \partial _x \phi \, \bigr) 
    - \bigl|\, \partial _x \tilde{U} \, \bigr|^{p-1} 
      \bigl( \, \partial _x \tilde{U} \, \bigr) \, 
   \Bigr)
   \, \biggr) \\
 &+(q-1)\int _0^{\phi} 
   \left( \lambda(\tilde{U}+\eta)-\lambda(\tilde{U}) \right)\left| 
   \, \eta \, \right|^{q-2}
   \, \mathrm{d}\eta \, \bigl( \, \partial _x \tilde{U} \, \bigr) \\ 
 &+\mu \, (q-1) \, \left|\, \phi \, \right|^{q-2} \partial _x \phi \, \\
 & \quad \times 
   \Bigl( \, 
   \bigl|\, \partial _x \tilde{U} + \partial _x \phi \, \bigr|^{p-1} 
   \bigl( \, \partial _x \tilde{U} + \partial _x \phi \, \bigr) 
   - \bigl|\, \partial _x \tilde{U} \, \bigr|^{p-1} 
     \bigl( \, \partial _x \tilde{U} \, \bigr) \, 
   \Bigr)\\
 &=\left|\, \phi \, \right|^{q-2} \phi \, F_{p}(U,U^r). 
 \end{aligned}
 \end{align}
Integrating (4.3) with respect to $x$, we have
 \begin{align}
 \begin{aligned}
 &\frac{1}{q} \frac{\mathrm{d}}{\mathrm{d}t} 
 \,\Vert \, \phi(t) \, \Vert_{L^q}^q \\
 &+\int ^{\infty }_{-\infty } (q-1) \int _0^{\phi} 
 \left( \lambda(\tilde{U}+\eta)-\lambda(\tilde{U}) \right)
 \left| \, \eta \, \right|^{q-2}
 \, \mathrm{d}\eta \, \bigl( \, \partial _x \tilde{U} \, \bigr) 
 \, \mathrm{d}x \\
 &+\mu \, (q-1) \, \int ^{\infty }_{-\infty }
   \left|\, \phi \, \right|^{q-2} \partial _x \phi \, \\
 & \quad \times 
   \Bigl( \, 
   \bigl|\, \partial _x \tilde{U} + \partial _x \phi \, \bigr|^{p-1} 
   \bigl( \, \partial _x \tilde{U} + \partial _x \phi \, \bigr) 
   - \bigl|\, \partial _x \tilde{U} \, \bigr|^{p-1} 
     \bigl( \, \partial _x \tilde{U} \, \bigr) \, 
   \Bigr) 
   \, \mathrm{d}x \\
 &=\int ^{\infty}_{-\infty}\left|\, \phi \, \right|^{q-2} \phi 
  \, F_{p}(U,U^r) \, \mathrm{d}x. 
 \end{aligned}
 \end{align} 
In order to estimate the second term on the left-hand side of 
(4.4), 
noting the shape of the flux function $f$, 
we divide the integral region of $x$ depending on the signs of
$\tilde U +\phi$, $\tilde U$ and $\phi$ as 
\begin{align*}
\begin{aligned}
&\int ^{\infty }_{-\infty } \biggl(\,  \int _0^{\phi} 
\left( \lambda \bigl(\tilde{U}+\eta \bigr)-\lambda \bigl(\tilde{U} \bigr) \right) 
\left| \, \eta \, \right|^{q-2}
\, \mathrm{d}\eta \biggr)
\, \bigl( \, \partial _x \tilde{U} \, \bigr) \, \mathrm{d}x 
\\
\qquad &= \int_{\tilde{U}+\phi \geq 0, \tilde{U} \geq 0, \phi \geq 0}
          +\int_{\tilde{U}+\phi \geq 0, \tilde{U} \geq 0, \phi \le 0}
          +\int_{\tilde{U}+\phi \geq 0, \tilde{U} < 0}
          + \int_{\tilde{U}+\phi < 0, \tilde{U} \geq 0}
\end{aligned}
\end{align*}
where we used the fact that the integral is clearly zero on the 
domain $\tilde{U}+\phi \le 0$ and $\tilde{U}\le 0$.
By Lagrange's mean-value theorem, we easily get as 
\begin{align*}
\begin{aligned}
&\left( \int ^{\infty }_{-\infty } \biggl(\,  \int _0^{\phi} 
\left( \lambda \bigl(\tilde{U}+\eta \bigr)-\lambda \bigl(\tilde{U} \bigr) \right) 
\left| \, \eta \, \right|^{q-2}
\, \mathrm{d}\eta \biggr)
\, \bigl( \, \partial _x \tilde{U} \, \bigr) \, \mathrm{d}x \right) (t)
\sim  G_{q}(t) 
\end{aligned}
\end{align*}
where $G_{q}=G_{q}(t)$ is defined in Proposition 4.1 
(cf. \cite{matsumura-yoshida}, \cite{yoshida}, \cite{yoshida'}). 
Next, by using the uniform boundedness, Lemma 3.1, 
and the following absolute equality 
with $p>1$, 
for any $a,b \in \mathbb{R}$, 
\begin{align*}
\begin{aligned}
&\left(\, |\, a\, |^{p-1}a - |\, b\, |^{p-1}b \, \right) 
 \left(\, a-b \, \right) \\
&= \frac{1}{2} \, 
   \left(\, |\, a\, |^{p-1} + |\, b\, |^{p-1} \, \right) 
   \left(\, a-b \, \right)^2 
   + \frac{1}{2} \, 
     \left(\, |\, a\, |^{p-1} - |\, b\, |^{p-1} \, \right) 
     \left(\, a^2-b^2 \, \right)\\
& \geq  \frac{1}{4} \, 
   \left(\, |\, a\, |^{p-1} + |\, b\, |^{p-1} + |\, a-b\, |^{p-1} \, \right) 
   \left(\, a-b \, \right)^2 
   + \frac{1}{2} \, 
     \left(\, |\, a\, |^{p-1} - |\, b\, |^{p-1} \, \right) 
     \left(\, a^2-b^2 \, \right),
\end{aligned}
\end{align*}
we have 
 \begin{align}
 \begin{aligned}
 &\frac{1}{q} \frac{\mathrm{d}}{\mathrm{d}t} 
 \,\Vert \, \phi(t) \, \Vert_{L^q}^q 
  +C_{p,q}^{-1} \, G_{q}(t) \\
 &+\frac{\mu \, (q-1)}{4} \, 
   \int _{-\infty}^{\infty} 
    | \, \phi \, |^{q-2} 
    \bigl( \, \partial _x \phi \, \bigr)^2 
    \left( \, 
    \bigl| \partial _x \phi \bigr|^{p-1} 
 + \bigl| \partial _x U \bigr|^{p-1} + \bigl| \partial _x U^r \bigr|^{p-1} \, 
   \right) \, \mathrm{d}x \\
 &+\frac{\mu \, (q-1)}{2} \, 
   \int ^{\infty }_{-\infty }
   | \, \phi \, |^{q-2} \, 
    \left| \, 
    \bigl| \partial _x \phi + \partial _x U^r \bigr|^{p-1} 
    - \bigl| \partial _x U^r \bigr|^{p-1}  
    \, \right| \\
 & \qquad \qquad \qquad \qquad \quad \; \; \, \, 
   \times 
    \left| \, 
    \bigl( \, \partial _x \phi + \partial _x U^r \, \bigr)^2 
    - \bigl( \, \partial _x U^r \, \bigr)^2 \, \right| 
    \, \mathrm{d}x \\
 &\leq \left| \, \int ^{\infty}_{-\infty}\left|\, \phi \, \right|^{q-2} \phi 
       \, F_{p}(U,U^r) \, \mathrm{d}x \, \right|. 
 \end{aligned}
 \end{align}
We note the right-hand side of (4.5) can be estimated as 
 \begin{align}
 \begin{aligned}
 &\left| \, 
  \int ^{\infty}_{-\infty} 
  \left|\, \phi \, \right|^{q-2} \phi 
  \, F_{p}(U,U^r) \, \mathrm{d}x \, 
  \right| 
 \leq \left| \, 
      \int ^{\infty}_{-\infty} 
      \left|\, \phi \, \right|^{q-1} \, 
      \widetilde{F_{p} }(U,U^r) 
      \, \mathrm{d}x \, 
      \right| \\
 & \qquad \; 
   +\int _{-\infty}^{\infty} 
    \left|\, \phi \, \right|^{q-2}
    \bigl| \, \partial _x \phi \, \bigr|
    \Bigl( \, 
    \bigl( \, \partial_x U + \partial_x U^{r} \, \bigr)^{p} 
    - \left( \, \partial_x U \, \right)^{p} \, 
      \Bigr) 
      \, \mathrm{d}x.
 \end{aligned}
 \end{align}
Thus, substituting (4.6) into (4.5), multiplying the inequality by 
$(1+t)^{\alpha}$ with $\alpha>0$ 
and integrating over $(0,t)$ with respect to the time, 
we complete the proof of Lemma 4.1.

\medskip

{\bf Proof of Proposition 4.1}.\
By using Lemma 4.1 and Lemma 4.2, 
we shall estimate the second term, the third term and the fourth term 
on the right-hand side of (4.2) as follows: 
for any $\epsilon>0$, 
\begin{align}
\begin{aligned}
 &\alpha \int ^t_0 (1+\tau )^{\alpha -1} 
  \| \,\phi(\tau) \,\|_{L^q}^q 
  \, \mathrm{d}\tau \\
 &\le C_{\alpha ,p,q} \int ^t_0 (1+\tau )^{\alpha -1} 
      \left( \, 
      \int _{-\infty}^{\infty} 
      | \, \phi \, |^{q-2} 
      \bigl| \, \partial _x \phi \, \bigr|^{p+1}
      \, \mathrm{d}x
      \right)^{\frac{q-2}{3p+q-1}} \\
 &\qquad \qquad \qquad \qquad \qquad \quad \times 
      \left( \, 
      \int _{-\infty}^{\infty} 
      | \, \phi \, |^{2} 
      \, \mathrm{d}x
      \, \right)^{\frac{pq+p+q-1}{3p+q-1}}
      \, \mathrm{d}\tau \\
 &\le \int ^t_0 
      \left( \, 
      (1+\tau )^{\alpha} 
      \int _{-\infty}^{\infty} 
      | \, \phi \, |^{q-2} 
      \bigl| \, \partial _x \phi \, \bigr|^{p+1}
      \, \mathrm{d}x
      \, \right)^{\frac{q-2}{3p+q-1}} \\
 &\qquad \quad \quad \quad \quad \times 
      C_{\alpha ,p,q}\,  (1+\tau )^{\alpha - 1 -\frac{\alpha (q-2)}{3p+q-1}}
         \| \,\phi(\tau) \,\|_{L^2}^{\frac{2(pq+p+q-1)}{3p+q-1}}
         \, \mathrm{d}\tau \\
 &\le \epsilon \int ^t_0 (1+\tau )^{\alpha} 
      \left( \, 
      \int _{-\infty}^{\infty} 
      | \, \phi \, |^{q-2} 
      \bigl| \, \partial _x \phi \, \bigr|^{p+1}
      \, \mathrm{d}x
      \right) \, \mathrm{d}\tau \\
 &\quad \quad \quad \quad \quad 
         +C_{\alpha ,p,q}(\epsilon ) 
         \int ^t_0 (1+\tau )^{\alpha - {\frac{3p+q-1}{3p+1}}}
         \| \,\phi(\tau ) \,\|_{L^2}^{\frac{2(pq+p+q-1)}{3p+1}} 
         \, \mathrm{d}\tau, 
\end{aligned}
\end{align}
 \begin{align}
 \begin{aligned}
 &q \int ^t_0 (1+\tau )^{\alpha} 
  \| \,\phi(\tau ) \,\|_{L^\infty}^{q-1} 
  \left|\left| \, \widetilde{F_{p} }(U,U^r)(\tau ) \, \right| \right|_{L^1}
  \, \mathrm{d}\tau \\
 &\le C_{p,q} \int ^t_0 (1+\tau )^{\alpha} 
      \left( \, 
      \int _{-\infty}^{\infty} 
      | \, \phi \, |^{q-1} 
      \bigl| \, \partial _x \phi \, \bigr|^{p+1}
      \, \mathrm{d}x
      \right)^{\frac{q-1}{3p+q-1}} \\
 &\qquad \quad \quad \quad \quad \times 
      \left( \, 
      \int _{-\infty}^{\infty} 
      | \, \phi \, |^{2} 
      \, \mathrm{d}x
      \, \right)^{\frac{p(q-1)}{3p+q-1}}
      \left|\left| \, \widetilde{F_{p} }(U,U^r)(\tau ) \, \right| \right|_{L^1}
      \, \mathrm{d}\tau \\
 &\le \int ^t_0 
      \left( \, 
      (1+\tau )^{\alpha} 
      \int _{-\infty}^{\infty} 
      | \, \phi \, |^{q-2} 
      \bigl| \, \partial _x \phi \, \bigr|^{p+1}
      \, \mathrm{d}x
      \, \right)^{\frac{q-2}{3p+q-1}} \\
 &\quad \: \: \: \times 
         C_{p,q}\,  (1+\tau )^{\alpha -\frac{\alpha (q-1)}{3p+q-1}}
         \| \,\phi(\tau) \,\|_{L^2}^{\frac{2p(q-1)}{3p+q-1}}
         \left|\left| \, \widetilde{F_{p} }(U,U^r)(\tau ) \, \right| \right|_{L^1}
         \, \mathrm{d}\tau \\
 &\le \epsilon \int ^t_0 (1+\tau )^{\alpha} 
      \left( \, 
      \int _{-\infty}^{\infty} 
      | \, \phi \, |^{q-1} 
      \bigl| \, \partial _x \phi \, \bigr|^{p+1}
      \, \mathrm{d}x
      \right) 
      \, \mathrm{d}\tau \\
 &\quad 
         + C_{p,q}(\epsilon) \int ^t_0 
          (1+\tau )^{\alpha} \| \,\phi(\tau ) \,\|
          _{L^2}^{{\frac{2}{3}}(q-1)}
          \left|\left| \, \widetilde{F_{p} }(U,U^r)(\tau ) \, \right| \right|_{L^1}
          ^{\frac{3p+q-1}{3p}} 
          \, \mathrm{d}\tau, 
 \end{aligned}
 \end{align}
  \begin{align}
 \begin{aligned}
 & \mu \, q \int ^t_0 (1+\tau )^\alpha 
   \| \,\phi(\tau ) \,\|_{L^\infty}^{p-2} \\
 & \quad \quad \, \times    
    \int _{-\infty}^{\infty} 
    \bigl| \, \partial _x \phi \, \bigr|
    \Bigl( \, 
    \bigl( \, \partial_x U + \partial_x U^{r} \, \bigr)^{p} 
    - \left( \, \partial_x U \, \right)^{p} \, 
     \Bigr) 
     \, \mathrm{d}x\mathrm{d}\tau \\
 &\le C_{p,q} \int ^t_0 (1+\tau )^{\alpha} 
      \left( \, 
      \int _{-\infty}^{\infty} 
      | \, \phi \, |^{q-1} 
      \bigl| \, \partial _x \phi \, \bigr|^{p+1}
      \, \mathrm{d}x
      \right)^{\frac{q-2}{3p+q-1}} \\
 &\quad \quad \, \times 
      \left( \, 
      \int _{-\infty}^{\infty} 
      | \, \phi \, |^{2} 
      \, \mathrm{d}x
      \, \right)^{\frac{p(q-2)}{3p+q-1}} \\
 &\quad \quad \, \times 
      \left( \, 
      \int _{-\infty}^{\infty} 
      \Bigl( \, 
      \bigl( \, \partial_x U + \partial_x U^{r} \, \bigr)^{p} 
      - \left( \, \partial_x U \, \right)^{p} \, 
      \Bigr)^{\frac{p+1}{p}} 
      \, \mathrm{d}x 
      \, \right)
      \mathrm{d}\tau \\
 &\le \int ^t_0 
      \left( \, 
      (1+\tau )^{\alpha} 
      \int _{-\infty}^{\infty} 
      | \, \phi \, |^{q-2} 
      \bigl| \, \partial _x \phi \, \bigr|^{p+1}
      \, \mathrm{d}x
      \, \right)^{\frac{q-2}{3p+q-1}} \\
 &\quad \: \: \: \times 
         C_{p,q}\,  (1+\tau )^{\alpha -\frac{\alpha (q-2)}{3p+q-1}}
         \| \,\phi(\tau) \,\|_{L^2}^{\frac{2p(q-2)}{3p+q-1}} \\
 &\quad \: \: \: \times 
         \left( \, 
      \int _{-\infty}^{\infty} 
      \Bigl( \, 
      \bigl( \, \partial_x U + \partial_x U^{r} \, \bigr)^{p} 
      - \left( \, \partial_x U \, \right)^{p} \, 
      \Bigr)^{\frac{p+1}{p}} 
      \, \mathrm{d}x 
      \, \right)
         \, \mathrm{d}\tau \\
 &\le \epsilon \int ^t_0 (1+\tau )^{\alpha} 
      \left( \, 
      \int _{-\infty}^{\infty} 
      | \, \phi \, |^{q-1} 
      \bigl| \, \partial _x \phi \, \bigr|^{p+1}
      \, \mathrm{d}x
      \right) 
      \, \mathrm{d}\tau \\
 &\quad \, 
         + C_{p,q}(\epsilon) \int ^t_0 
          (1+\tau )^{\alpha} \| \,\phi(\tau ) \,\|
          _{L^2}^{\frac{2p(q-2)}{3p+1}}\\
 &\quad \, \, \: \: \: \times 
      \left( \, 
      \int _{-\infty}^{\infty} 
      \Bigl( \, 
      \bigl( \, \partial_x U + \partial_x U^{r} \, \bigr)^{p} 
      - \left( \, \partial_x U \, \right)^{p} \, 
      \Bigr)^{\frac{p+1}{p}} 
      \, \mathrm{d}x 
      \, \right)^{\frac{3p+q-1}{3p+1}}
          \, \mathrm{d}\tau. 
 \end{aligned}
 \end{align}
Substituting (4.7), (4.8) and (4.9) into (4.2), 
we have 
\begin{align}
\begin{aligned}
 &(1+t)^\alpha  \| \,\phi(t) \,\|_{L^q}^q 
  + \int ^t_0 (1+\tau )^\alpha  G_q(\tau ) \, \mathrm{d}\tau \\
 & \quad 
   + \int ^t_0 (1+\tau )^\alpha  \int _{-\infty}^{\infty} 
    | \, \phi \, |^{q-2} 
    \bigl( \, \partial _x \phi \, \bigr)^2 \\
 & \qquad \qquad \qquad \qquad \quad \quad \; \, \, 
    \times 
    \left( 
    \bigl| \partial _x \phi \bigr|^{p-1} 
 + \bigl| \partial _x U \bigr|^{p-1} + \bigl| \partial _x U^r \bigr|^{p-1}  
   \right) \, \mathrm{d}x \mathrm{d}\tau  \\
 & \quad 
  + \int ^t_0 (1+\tau )^\alpha  \int _{-\infty}^{\infty} 
    | \, \phi \, |^{q-2} \, 
    \left| \, 
    \bigl| \partial _x \phi + \partial _x \tilde{U} \bigr|^{p-1} 
    - \bigl| \partial _x \tilde{U} \bigr|^{p-1}  
    \, \right| \\
 & \qquad \qquad \qquad \quad \quad \quad \; \, \, 
   \times 
    \left| \, 
    \bigl( \, \partial _x \phi + \partial _x \tilde{U} \, \bigr)^2 
    - \bigl( \, \partial _x \tilde{U} \, \bigr)^2 \, \right| 
    \, \mathrm{d}x \mathrm{d}\tau  \\
&\leq C_{\alpha,p,q} \| \,\phi_0 \,\|_{L^p}^p 
 + C_{\alpha ,p,q}(\epsilon ) 
         \int ^t_0 (1+\tau )^{\alpha - {\frac{3p+q-1}{3p+1}}}
         \| \,\phi(\tau ) \,\|_{L^2}^{\frac{2(pq+p+q-1)}{3p+1}} 
         \, \mathrm{d}\tau \\
 &\qquad \quad \; \; \: \,  
 + C_{p,q}(\epsilon) \int ^t_0 
          (1+\tau )^{\alpha} \| \,\phi(\tau ) \,\|
          _{L^2}^{{\frac{2}{3}}(q-1)}
          \left|\left| \, \widetilde{F_{p} }(U,U^r)(\tau ) \, \right| \right|_{L^1}
          ^{\frac{3p+q-1}{3p}} 
          \, \mathrm{d}\tau \\
 &\qquad \quad \; \; \: \,  
 + C_{p,q}(\epsilon) \int ^t_0 
          (1+\tau )^{\alpha} \| \,\phi(\tau ) \,\|
          _{L^2}^{\frac{2p(q-2)}{3p+1}}\\ 
 &\qquad \quad \quad \, \, \: \: \: \; \; \; \; \; \times 
      \left( \, 
      \int _{-\infty}^{\infty} 
      \Bigl( \, 
      \bigl( \, \partial_x U + \partial_x U^{r} \, \bigr)^{p} 
      - \left( \, \partial_x U \, \right)^{p} \, 
      \Bigr)^{\frac{p+1}{p}} 
      \, \mathrm{d}x 
      \, \right)^{\frac{3p+q-1}{3p+1}}
          \, \mathrm{d}\tau. 
 \end{aligned}
 \end{align}
By using the $L^2$-boundedness of $\phi$, Lemma 3.2, 
we first get 
 \begin{align}
 \begin{aligned}
 \| \, \phi (t)\, \| _{L^2}^2 
 \leq C_{p}(\phi_0). 
 \end{aligned}
 \end{align}
By using Lemma 2.2, Lemma 2.3 and Lemma 3.3, 
we also get 

\medskip

\noindent
{\bf Lemma 4.3.}\quad {\it
For any fixed $p \in (1, \infty)$, 
we have the following time-decay estimates. 

\noindent
{\rm (1)}\ \ For any $\delta \in (0, 1)$, there exists positive 
constants $C_{p}$, $C_{\delta }$ and $T_{0}$ such that 
\begin{align*}
\begin{aligned}
\left|\left| \, \widetilde{F_{p} }(U,U^r)(t) \, \right| \right|_{L^1} 
\leq C_{p} \left( \, (1+t)^{-{\frac{2p}{p+1}}} 
     + C_{\delta }(1+t)^{-2(1-\delta )} \, \right) 
     \quad \bigl( \, t \ge T_{0} \, \bigr).
\end{aligned}
\end{align*}

\noindent
{\rm (2)}\ \ There exists a positive constant $C_{p}$ 
such that 
\begin{align*}
\begin{aligned}
\left( \, 
\int _{-\infty}^{\infty} 
      \Bigl( \, 
      \bigl( \, \partial_x U + \partial_x U^{r} \, \bigr)^{p} 
      - \left( \, \partial_x U \, \right)^{p} \, 
      \Bigr)^{\frac{p+1}{p}} 
      \, \mathrm{d}x 
\, \right) (t)
\leq C_{p}(1+t)^{-1}
      \quad \bigl( t \ge 0 \bigr).
\end{aligned}
\end{align*}
}

\medskip

\noindent
We estimate the each terms on the right-hand side of (4.10) 
as follows: 
\begin{align}
\begin{aligned}
&C_{\alpha ,p,q} \, 
         \int ^t_0 (1+\tau )^{\alpha - {\frac{3p+q-1}{3p+1}}} 
         \| \,\phi(\tau ) \,\|_{L^2}^{\frac{2(pq+p+q-1)}{3p+1}} 
         \, \mathrm{d}\tau \\
&\leq C_{\alpha ,p,q} 
      \bigl( \, C_{p}(\phi_0) \, \bigr)^{\frac{pq+p+q-1}{3p+1}} 
      \int ^t_0 (1+\tau )^{\alpha - {\frac{3p+q-1}{3p+1}}} 
      \, \mathrm{d}\tau \\
&\leq C_{\alpha ,p,q} 
      \bigl( \, C_{p}(\phi_0) \, \bigr)^{\frac{pq+p+q-1}{3p+1}} 
      (1+t )^{\alpha - {\frac{q-2}{3p+1}}}, 
\end{aligned}
\end{align}
\begin{align}
\begin{aligned}
&C_{p,q} \int ^t_0 
          (1+\tau )^{\alpha} \| \,\phi(\tau ) \,\|
          _{L^2}^{{\frac{2}{3}}(q-1)}
          \left|\left| \, \widetilde{F_{p} }(U,U^r)(\tau ) \, \right| \right|_{L^1}
          ^{\frac{3p+q-1}{3p}} 
          \, \mathrm{d}\tau \\
&\leq C_{p,q} 
      \bigl( \, C_{p}(\phi_0) \, \bigr)^{{\frac{1}{3}}(q-1)} 
      \int ^t_0 (1+\tau )^{\alpha - {\frac{2(3p+q-1)}{3(p+1)}}} 
      \, \mathrm{d}\tau \\
&\leq C_{p,q} 
      \bigl( \, C_{p}(\phi_0) \, \bigr)^{{\frac{1}{3}}(q-1)} 
      (1+t )^{\alpha - {\frac{3p+2q-5}{3(p+1)}}}, 
\end{aligned}
\end{align}
\begin{align}
\begin{aligned}
&C_{p,q} \int ^t_0 
          (1+\tau )^{\alpha} \| \,\phi(\tau ) \,\|
          _{L^2}^{\frac{2p(q-2)}{3p+1}}\\ 
&\quad \quad \, \, \: \: \: \; \; \; \; \; \times 
      \left( \, 
      \int _{-\infty}^{\infty} 
      \Bigl( \, 
      \bigl( \, \partial_x U + \partial_x U^{r} \, \bigr)^{p} 
      - \left( \, \partial_x U \, \right)^{p} \, 
      \Bigr)^{\frac{p+1}{p}} 
      \, \mathrm{d}x 
      \, \right)^{\frac{3p+q-1}{3p+1}} (\tau )
      \, \mathrm{d}\tau \\
&\leq C_{p,q} 
      \bigl( \, C_{p}(\phi_0) \, \bigr)^{{\frac{p(q-2)}{3p+1}}} 
      \int ^t_0 (1+\tau )^{\alpha - {\frac{3p+q-1}{3p+1}}} 
      \, \mathrm{d}\tau \\
&\leq C_{p,q} 
      \bigl( \, C_{p}(\phi_0) \, \bigr)^{{\frac{p(q-2)}{3p+1}}} 
      (1+t )^{\alpha - {\frac{q-2}{3p+1}}}. 
\end{aligned}
\end{align}
Substituting (4.12), (4.13) and (4.14) into (4.10), 
we get (4.1). 
Thus the proof of Proposition 4.1 is complete. 
In particular, it follows that 
\begin{align}
\begin{aligned}
\| \,\phi(t) \,\|_{L^q} 
\leq C( \, p, \,q, \, \phi_0 \, ) \, 
     (1+t)^{-{\frac{1}{3p+1}} \left( 1 - {\frac{2}{q}} \right)}
\end{aligned}
\end{align}
for $2 \leq q < \infty$. 

\medskip

{\bf Proof of Theorem 3.1.}\ 
We already have proved the decay estimate 
of $\| \,\phi(t) \,\|_{L^q}$ with $2 \le q < \infty$. 
Therefore we only 
show the $L^{\infty}$-estimate. 
We first note by Lemma 2.2 and Lemma 2.3 that 
\begin{align}                
\begin{aligned}
&\bigl|\bigl| 
 \, \partial _x \phi (t)\, 
 \bigr|\bigr| _{L^{p+1}}^{p+1}\\
&\leq \bigl|\bigl| 
      \, \partial _x u (t)\, 
      \bigr|\bigr| _{L^{p+1}}^{p+1} 
      + \bigl|\bigl| 
        \, \partial _x U^r (t)\, 
        \bigr|\bigr| _{L^{p+1}}^{p+1} \\
&\leq C\bigl( \, \| \, \phi_0 \, \|_{L^2}, \: 
       \| \, \partial _x u_0 \, \|_{L^{p+1}} 
       \, \bigr) 
       + C_{p}(1+t)^{-1}. 
\end{aligned}
\end{align}
We use the following Gagliardo-Nirenberg inequality: 
\begin{align}
\| \,\phi(t) \,\|_{L^\infty }
\le C_{q,\theta}
\| \,\phi(t) \,\|_{L^q}^{1-\theta}
\| \,\partial _x \phi(t) \,\|_{L^{p+1}}^{\theta}
\end{align}
for any 
$(q,\theta)\in [\, 1,\infty)\times (0,1\, ]$ satisfying 
$$
\frac{p}{p+1} \, {\theta}=(1-\theta) \, \frac{1}{q}. 
$$
Substituting (4.15) and (4.16) into (4.17), 
we have 
\begin{align}
\begin{aligned}
\| \,\phi(t) \,\|_{L^\infty }
&\le C( \, p, \,q, \, \theta, \, \phi_0, \, \partial _x u_0 \, ) \, 
     (1+t)^{-{\frac{1}{3p+1}} \left( 1 - {\frac{2}{q}} \right)(1-\theta)}\\
&\le C( \, p, \, \theta, \, \phi_0, \, \partial _x u_0 \, ) \, 
     (1+t)^{-{\frac{1}{3p+1}} + {\frac{\theta}{p+1}}}
\end{aligned}
\end{align}
for $\theta \in (0,1\, ]$. 
Consequently, we do complete the proof of Theorem 3.1.
\bigskip 

\noindent
\section{Time-decay estimates with $1 \le q \le \infty$ }
In this section, 
we show the time-decay estimates 
with $1 \le q \le \infty$ and 
time-decay estimate for the higher order derivative 
in the $L^{p+1}$-norm, 
in the case where $\phi_0 \in L^1 \cap L^2$ 
with $\partial _x u_0 \in L^{p+1}$, 
that is, Theorem 3.2. 
Then, we 
first establish the $L^1$-estimate 
to the solution $\phi$ of the reformulated Cauchy problem (3.6). 
To do that, we use the Friedrichs mollifier $\rho_\delta \ast $, 
where, 
$\rho_\delta(\phi):=\frac{1}{\delta}\rho \left( \frac{\phi}{\delta}\right)$
with 
\begin{align*}
\begin{aligned}
&\rho \in C^{\infty}_0(\mathbb{R}), \quad
\rho (\phi)\geq 0 \quad (\phi \in \mathbb{R}),
\\
&\mathrm{supp} \{\rho \} \subset 
\left\{\phi \in \mathbb{R}\: \left|\:  |\,\phi \, |\le 1 \right. \right\},\quad  
\int ^{\infty}_{-\infty} \rho (\phi)\, \mathrm{d}\phi=1.
\end{aligned}
\end{align*}

Some useful properties of the mollifier are as follows. 
\medskip

\noindent
{\bf Lemma 5.1.}\quad{\it

\noindent
{\rm (i)}\ 
$\displaystyle{\lim_{\delta \to 0}\, 
\left( \rho_\delta \ast \mathrm{sgn} \right)(\phi)}
= \mathrm{sgn} (\phi)
\qquad (\phi \in \mathbb{R}),$

\noindent
{\rm (ii)}\ 
$\displaystyle{\lim_{\delta \to 0}\, 
\int^{\phi}_0\left( \rho_\delta \ast \mathrm{sgn} \right)(\eta)
\, \mathrm{d}\eta}
=|\, \phi \, |
\qquad (\phi \in \mathbb{R}),$

\medskip

\noindent
{\rm (iii)}\ 
$\Big. \left( \rho_\delta \ast \mathrm{sgn} \right) \Bigr|_{\phi=0} =0,$

\medskip

\noindent
{\rm (iv)}\ 
$\displaystyle{ \frac{\mathrm{d}}{\mathrm{d}\phi}\, 
\left( \rho_\delta \ast \mathrm{sgn} \right)(\phi)}
=2\, \rho_\delta(\phi)
\ge 0
\qquad (\phi \in \mathbb{R}),$

\medskip

\noindent
where 
$$
\left( \rho_\delta \ast \mathrm{sgn} \right)(\phi)
:=\int^{\infty}_{-\infty}
\rho_\delta(\phi-y)\, \mathrm{sgn}(y)\, \mathrm{d}y
\qquad (\phi \in \mathbb{R})
$$
and 
\begin{equation*}
\mathrm{sgn}(\phi):=\left\{
\begin{array}{ll}
-1& \quad \; \bigl(\, \phi < 0 \, \bigr),\\[3pt]
\, \, 0& \quad \; \bigl(\, \phi = 0 \, \bigr),\\[3pt]
\, \, 1& \quad \; \bigl(\, \phi > 0 \, \bigr).
\end{array}
\right.
\end{equation*} 
}

\noindent
Making use of Lemma 5.1, 
we obtain the following $L^1$-estimate. 

\medskip

\noindent
{\bf Proposition 5.1.}\quad {\it
Assume that the same assumptions in Theorem 3.2. 
For any $p>1$ and any $\epsilon>0$, 
there exist positive constants $C_{p}$ and $C_{\epsilon}$ 
such that 
the unique global solution in time $\phi$ 
of the Cauchy problem {\rm(3.6)} 
satisfies the following $L^1$-estimate 
\begin{align}
\begin{aligned}
\| \, \phi(t) \, \|_{L^1}
\le \| \,\phi_0 \,\|_{L^1} + C_{p}(1+t)^{-\frac{2p}{p+1}} 
    + C_{\epsilon}(1+t)^{-2(1-\epsilon)} 
\quad \bigl( \, t \ge T_0 \, \bigr) 
\end{aligned}
\end{align}
for any $\epsilon>0$.
}

\medskip

\noindent
{\bf Proof of Proposition 5.1.}\quad
Multiplying the equation in the problem (3.6) by 
$\left( \rho_\delta \ast \mathrm{sgn} \right)(\phi)$, 
we obtain the divergence form 
\medskip
\begin{align}
\begin{aligned}
&\partial_t\left( 
\int^{\phi}_0 
\left( \rho_\delta \ast \mathrm{sgn} \right)(\eta)
\, \mathrm{d}\eta 
\right) \\
&+\partial _x 
  \biggl( \, 
  \left( \rho_\delta \ast \mathrm{sgn} \right)(\phi) \, 
  \bigl( f(\tilde{U}+\phi)-f(\tilde{U}) \bigr) 
  \biggr) \\ 
&+\partial _x \left( 
  -\int _0^{\phi} 
  \bigl( f(\tilde{U}+\eta)-f(\tilde{U}) \bigr)
  \frac{\mathrm{d}\left( \rho_\delta \ast \mathrm{sgn} \right)}
  {\mathrm{d}\phi}
  (\eta)
  \, \mathrm{d}\eta 
  \, \right) \\
&+\partial _x \biggl( \, 
  -\mu \, \left( \rho_\delta \ast \mathrm{sgn} \right)(\phi) \biggr. \\
& \quad \times \biggl. 
   \Bigl( \, 
    \bigl|\, \partial _x \tilde{U} + \partial _x \phi \, \bigr|^{p-1} 
    \bigl( \, \partial _x \tilde{U} + \partial _x \phi \, \bigr) 
    - \bigl|\, \partial _x \tilde{U} \, \bigr|^{p-1} 
      \bigl( \, \partial _x \tilde{U} \, \bigr) \, 
   \Bigr)
  \, \biggr) \\
&+\int _0^{\phi} 
  \bigl( \lambda(\tilde{U}+\eta)-\lambda(\tilde{U}) \bigr) 
  \frac{\mathrm{d}\left( \rho_\delta \ast \mathrm{sgn} \right)}
  {\mathrm{d}\phi}(\eta) 
  \, \mathrm{d}\eta \, \bigl( \, \partial _x \tilde{U} \, \bigr) \\ 
&+\mu \, 
  \frac{\mathrm{d}\left( \rho_\delta \ast \mathrm{sgn} \right)}
  {\mathrm{d}\phi}(\phi) \, 
  \partial _x \phi \\
& \quad \times 
  \Bigl( \, 
    \bigl|\, \partial _x \tilde{U} + \partial _x \phi \, \bigr|^{p-1} 
    \bigl( \, \partial _x \tilde{U} + \partial _x \phi \, \bigr) 
    - \bigl|\, \partial _x \tilde{U} \, \bigr|^{p-1} 
      \bigl( \, \partial _x \tilde{U} \, \bigr) \, 
   \Bigr) \\
& = \left( \rho_\delta \ast \mathrm{sgn} \right)(\phi) \, 
   F_{p}(U,U^r). 
\end{aligned}
\end{align}
Integrating (5.2) with respect to $x$ and $t$, we have
\begin{align}
\begin{aligned}
&\int ^{\infty}_{-\infty} \int^{\phi (t) }_0 
 \left( \rho_\delta \ast \mathrm{sgn} \right)(\eta)
 \, \mathrm{d}\eta \, \mathrm{d}x \\
&+ \int _0^t \int ^{\infty}_{-\infty} \int^{\phi}_0
   \bigl( \lambda(\tilde{U}+\eta)-\lambda(\tilde{U}) \bigr) 
   \frac{\mathrm{d}\left( \rho_\delta \ast \mathrm{sgn} \right)}
   {\mathrm{d}\phi}(\eta)
   \, \mathrm{d}\eta \, \bigl( \, \partial _x \tilde{U} \, \bigr)
   \, \mathrm{d}x \mathrm{d}\tau \\
&+ \frac{\mu}{2} \int _0^t \int ^{\infty}_{-\infty} 
   \frac{\mathrm{d}\left( \rho_\delta \ast \mathrm{sgn} \right)}
   {\mathrm{d}\phi}(\phi) \, 
   \bigl( \, \partial _x \phi \, \bigr)^2 
    \left( \, 
    \bigl| \partial _x \phi \bigr|^{p-1} 
    + \bigl| \partial _x \tilde{U} \bigr|^{p-1}  
    \, \right)
   \, \mathrm{d}x \mathrm{d}\tau \\
&+ \frac{\mu}{2} \int _0^t \int ^{\infty}_{-\infty} 
   \frac{\mathrm{d}\left( \rho_\delta \ast \mathrm{sgn} \right)}
   {\mathrm{d}\phi}(\phi) \, 
   \left| \, 
    \bigl| \partial _x \phi + \partial _x \tilde{U} \bigr|^{p-1} 
    - \bigl| \partial _x \tilde{U} \bigr|^{p-1}  
    \, \right| \\
& \qquad \qquad \qquad \qquad \quad \quad \quad \; \, \, 
   \times 
    \left| \, 
    \bigl( \, \partial _x \phi + \partial _x \tilde{U} \, \bigr)^2 
    - \bigl( \, \partial _x \tilde{U} \, \bigr)^2 \, \right| 
   \, \mathrm{d}x \mathrm{d}\tau \\
&= \int ^{\infty}_{-\infty} \int^{\phi _0}_0 
   \left( \rho_\delta \ast \mathrm{sgn} \right)(\eta)
   \, \mathrm{d}\phi \, \mathrm{d}x \\
&\qquad \qquad  \qquad \; 
 + \int _0^t \int ^{\infty}_{-\infty} 
   \left( \rho_\delta \ast \mathrm{sgn} \right)(\phi) \, 
   F_{p}(U,U^r)
   \, \mathrm{d}x \mathrm{d}\tau .
\end{aligned}
\end{align}
By using Lemma 5.1, we first note that for $t\in [\, 0,\infty)$,
\begin{align}
\left| \, 
\int^{\phi (t)}_0\left( \rho_\delta \ast \mathrm{sgn} \right)(\eta)
\, \mathrm{d}\eta 
\, \right|
\le \left( \rho_\delta \ast \mathrm{sgn} \right)
    (|\,  \phi (t) \,|)\, |\,  \phi (t) \,|
\le |\,  \phi (t) \,|, 
\end{align}
\begin{align}
\displaystyle {\lim_{\delta\rightarrow 0}}
\int ^{\infty}_{-\infty} \int^{\phi (t) }_0
\left( \rho_\delta \ast \mathrm{sgn} \right)(\eta)
\, \mathrm{d}\eta 
= \| \,\phi(t) \,\|_{L^1}. 
\end{align}
We also note the following (the proof is similar to the one in \cite{yoshida}). 
\medskip

\noindent
{\bf Lemma 5.2.}\quad{\it
It holds that 
\begin{align}
\begin{aligned}
& \int ^{\infty}_{-\infty} \int^{\phi 
(t) }_0
  \left( \lambda(\tilde{U}+\eta)-\lambda(\tilde{U}) \right) 
  \frac{\mathrm{d}\left( \rho_\delta \ast \mathrm{sgn} \right)}
  {\mathrm{d}\phi}(\eta)
  \, \mathrm{d}\eta \, \bigl( \, \partial _x \tilde{U} \, \bigr)
  \, \mathrm{d}x \\
& \geq C^{-1} \left( \, \int_{\tilde{U}+\phi \geq 0, \tilde{U} \geq 0 } 
              \left|\, \int_{0}^{| \phi |} 
              \eta \, \rho _{\delta} (\eta) \, \mathrm{d}\eta \, \right|
              \bigl( \partial _x \tilde{U} \bigr) 
              \, \mathrm{d}x \, \right)(t)\\
& \quad        
       + C^{-1} \left( \, \int_{\tilde{U}+\phi \geq 0, \tilde{U} < 0 } 
                \left|\, \int_{0}^{| \phi |} 
                \bigl(\tilde{U}+\eta \bigr) \, \rho _{\delta} (\eta) \, 
                \mathrm{d}\eta \, \right|
                \bigl( \partial _x \tilde{U} \bigr) 
                \, \mathrm{d}x \, \right)(t)\\
& \quad 
       + C^{-1} \left( \, \int_{\tilde{U}+\phi < 0, \tilde{U} \geq 0 } 
                \left|\, \int_{0}^{| \phi |} 
                \tilde{U} \, \rho _{\delta} (\eta) \, 
                \mathrm{d}\eta \, \right|
                \bigl( \partial _x \tilde{U} \bigr) 
                \, \mathrm{d}x  \, \right)(t)
                \geq 0
                \qquad \bigl( t \ge 0 \bigr). 
\end{aligned}
\end{align}
}
\medskip

\noindent
So we can get 
\begin{align}
\begin{aligned}
&\| \,\phi (t) \,\|_{L^1} 
\le \| \,\phi_0 \,\|_{L^1} 
     + \displaystyle {\lim_{\delta\rightarrow 0}}
     \int^t_0 \left|\, \int ^{\infty}_{-\infty} 
     \left( \rho_\delta \ast \mathrm{sgn} \right)(\phi) \, 
     F_{p}(U,U^r)
   \, \mathrm{d}x \, \right|\, \mathrm{d}\tau \\ 
& \qquad \; \; \; \; \; \; \; \, 
     \le \| \,\phi_0 \,\|_{L^1} 
     + \displaystyle {\lim_{\delta\rightarrow 0}}
     \int^t_0 \left|\, \int ^{\infty}_{-\infty} 
     \left( \rho_\delta \ast \mathrm{sgn} \right)(\phi) \, 
     \widetilde{F_{p}}(U,U^r)
     \, \mathrm{d}x \, \right|\, \mathrm{d}\tau \\
& \qquad \qquad \quad \; \, 
   + \mu \, \displaystyle {\lim_{\delta\rightarrow 0}}
     \int^t_0 \biggl|\, \int ^{\infty}_{-\infty} 
     \left( \rho_\delta \ast \mathrm{sgn} \right)(\phi) \, \biggr.\\
& \quad \; 
     \biggl. \, \times 
     \partial _x \Bigl( \, 
     \bigl|\, \partial _x U + \partial _x U^r \, \bigr|^{p-1} 
     \bigl( \, \partial _x U + \partial _x U^r \, \bigr) 
     - \bigl|\, \partial _x U \, \bigr|^{p-1} 
      \bigl( \, \partial _x U \, \bigr) \, 
     \Bigr)
     \, \mathrm{d}x \, \biggr|\, \mathrm{d}\tau 
\end{aligned}
\end{align}
Noting by the asymptotic prorerties of $X(t)\, (t \ge T_0)$, Lemma 3.3, that 
\begin{align*}
\begin{aligned}
&\left( \, \int ^{X(t)}_{-\infty} + \int ^{\infty}_{X(t)} \, \right) 
 \bigl| \, \lambda(U+U^r)-\lambda(U^r) \, \bigr| \, \partial _x U^r
 \, \mathrm{d}x \\
&\qquad \qquad \qquad \quad \; \; \; \: \, 
\leq C_{p}(1+t)^{-1-\frac{p-1}{p+1}} + C_{\epsilon}(1+t)^{-1-(1-2\epsilon)} 
\quad \left( \, \epsilon \in (0,1) \, \right)
\end{aligned}
\end{align*}
and
$$
\int ^{\infty}_{-\infty} \lambda(U+U^r) \, \partial _x U
\, \mathrm{d}x 
\leq C_{p}(1+t)^{-1-\frac{p-1}{p+1}} + C_{\epsilon}(1+t)^{-1-(1-2\epsilon)} 
\quad \left( \, \epsilon \in (0,1) \, \right), 
$$
we immediately get 
\begin{align}
\begin{aligned}
& \displaystyle {\lim_{\delta\rightarrow 0}}\, 
  \left|\, \int ^{\infty}_{-\infty} 
  \left( \rho_\delta \ast \mathrm{sgn} \right)(\phi) \, 
  \widetilde{F_{p}}(U,U^r)
  \, \mathrm{d}x \, \right|\, (t) \\
& \leq  \left(\, \int ^{\infty}_{-\infty} 
  \bigl| \, \mathrm{sgn} (\phi) \, \bigr| \, 
  \Bigl| \, 
  \widetilde{F_{p}}(U,U^r)
  \, \Bigr| 
  \, \mathrm{d}x \, \right)\, (t) \\
& \leq C_{p}(1+t)^{-\frac{2p}{p+1}} 
    + C_{\epsilon}(1+t)^{-2(1-\epsilon)} 
  \quad \bigl( \, t \ge T_0 \, \bigr). 
\end{aligned}
\end{align}
Next, we estimate 
\begin{align*}
\begin{aligned}
&\biggl|\, \int ^{\infty}_{-\infty} 
 \left( \rho_\delta \ast \mathrm{sgn} \right)(\phi) \, 
 \partial _x \Bigl( \, 
 \bigl|\, \partial _x U + \partial _x U^r \, \bigr|^{p-1} 
 \bigl( \, \partial _x U + \partial _x U^r \, \bigr) 
 - \bigl|\, \partial _x U \, \bigr|^{p-1} 
 \bigl( \, \partial _x U \, \bigr) \, 
 \Bigr)
 \, \mathrm{d}x \, \biggr|\\
& \leq p \, \biggl|\, \int ^{\infty}_{-\infty} 
       \left( \rho_\delta \ast \mathrm{sgn} \right)(\phi) \, 
       \Bigl( \, 
       \bigl(\, \partial _x U + \partial _x U^r \, \bigr)^{p-1} 
       - \bigl(\, \partial _x U \, \bigr)^{p-1} \, 
       \Bigr) \, \partial _x^2 U
       \, \mathrm{d}x \, \biggr| \\
& \quad \, 
  + p \, \biggl|\, \int ^{\infty}_{-\infty} 
       \left( \rho_\delta \ast \mathrm{sgn} \right)(\phi) \, 
       \bigl(\, \partial _x U + \partial _x U^r \, \bigr)^{p-1} 
       \, \partial _x^2 U^r
       \, \mathrm{d}x \, \biggr|.
\end{aligned}
\end{align*}
By using Lagrange's mean-value theorem, we have
\begin{align*}
\begin{aligned}
& \biggl|\, \int ^{\infty}_{-\infty} 
  \left( \rho_\delta \ast \mathrm{sgn} \right)(\phi) \, 
  \Bigl( \, 
  \bigl(\, \partial _x U + \partial _x U^r \, \bigr)^{p-1} 
  - \bigl(\, \partial _x U \, \bigr)^{p-1} \, 
  \Bigr) \, \partial _x^2 U
  \, \mathrm{d}x \, \biggr|\\
& \leq \left\{\begin{array} {ll} 
       \displaystyle{
       (p-1) \, \biggl|\, \int ^{\infty}_{-\infty} 
       \left( \rho_\delta \ast \mathrm{sgn} \right)(\phi) \, 
       \bigl(\, \partial _x U \, \bigr)^{p-2} \, 
       \partial _x U^r \, \partial _x^2 U 
       \, \mathrm{d}x \, \biggr| 
       \quad \quad \quad \; \bigl( \, 1<p<2 \, \bigr)
       }, \\[15pt]
       \displaystyle{
       (p-1) \, \biggl|\, \int ^{\infty}_{-\infty} 
       \left( \rho_\delta \ast \mathrm{sgn} \right)(\phi) \, 
       \bigl(\, \partial _x U + \partial _x U^r \, \bigr)^{p-2} \, 
       \partial _x U^r \, \partial _x^2 U 
       \, \mathrm{d}x \, \biggr|
       \quad \; \; \bigl( \, p \geq 2 \, \bigr). 
       }
       \end{array}
       \right.\, 
\end{aligned}
\end{align*}
By using Lemma 2.2 and Lemma 2.3, we have for $1<p<2$, 
\begin{align}
\begin{aligned}
& \displaystyle {\lim_{\delta\rightarrow 0}}\, 
  \biggl|\, \int ^{\infty}_{-\infty} 
  \left( \rho_\delta \ast \mathrm{sgn} \right)(\phi) \, 
  \bigl(\, \partial _x U  \, \bigr)^{p-2} \, 
  \partial _x U^r \, \partial _x^2 U 
  \, \mathrm{d}x \, \biggr|\, (t) \\
& \leq \lint ^{\sqrt{\frac{A}{B}}{(1+t)}^{\frac{1}{p+1}} }_{-\sqrt{\frac{A}{B}}{(1+t)}^{\frac{1}{p+1}} }       
       \bigl(\, \partial _x U  \, \bigr)^{p-2} \, 
       \partial _x U^r \, | \, \partial _x^2 U \,|
       \, \mathrm{d}x \\
& \leq C_{p} (1+t)^{-\left( 1+\frac{p}{p+1} \right)}
       \lint ^{\sqrt{\frac{A}{B}}{(1+t)}^{\frac{1}{p+1}} }_{0} 
       \frac{x}{(1+t)^{\frac{1}{p+1}}} 
       \, \mathrm{d}x \\
& \leq C_{p} (1+t)^{-\frac{2p}{p+1}}
  \quad \bigl( \, t \ge T_0 \, \bigr), 
\end{aligned}
\end{align}
and for $p \geq 2$, 
\begin{align}
\begin{aligned}
& \displaystyle {\lim_{\delta\rightarrow 0}}\, 
  \biggl|\, \int ^{\infty}_{-\infty} 
  \left( \rho_\delta \ast \mathrm{sgn} \right)(\phi) \, 
  \bigl(\, \partial _x U + \partial _x U^r \, \bigr)^{p-2} \, 
  \partial _x U^r \, \partial _x^2 U 
  \, \mathrm{d}x \, \biggr|\, (t) \\
& \leq C_{p} \int ^{\infty}_{-\infty} 
       \bigl(\, \partial _x U  \, \bigr)^{p-2} \, 
       \partial _x U^r \, | \, \partial _x^2 U \,|
       \, \mathrm{d}x
       + C_{p} \int ^{\infty}_{-\infty} 
         \bigl(\, \partial _x U^r  \, \bigr)^{p-1} \, 
         | \, \partial _x^2 U \,| 
         \, \mathrm{d}x \\
& \leq C_{p} (1+t)^{-\frac{2p}{p+1}} 
       + C_{p} (1+t)^{-\left( p-1+\frac{2}{p+1} \right)} \\
& \quad \, \times 
         \lint ^{\sqrt{\frac{A}{B}}{(1+t)}^{\frac{1}{p+1}} }_{-\sqrt{\frac{A}{B}}{(1+t)}^{\frac{1}{p+1}} }       
         \left( \, 
         A-B \left(\frac{x}{(1+t)^{\frac{1}{p+1}}} \right)^2 \, 
         \right)^{-{\frac{p-2}{p-1}}} 
         \frac{| \, x \, |}{(1+t)^{\frac{1}{p+1}}}
         \, \mathrm{d}x \\
& \leq C_{p} (1+t)^{-\frac{2p}{p+1}} 
       + C_{p} (1+t)^{-\frac{p^2}{p+1}} \\
& \leq C_{p} (1+t)^{-\frac{2p}{p+1}}
  \quad \bigl( \, t \ge T_0 \, \bigr). 
\end{aligned}
\end{align}
Similarly, we have 
\begin{align}
\begin{aligned}
& \displaystyle {\lim_{\delta\rightarrow 0}}\, 
  \biggl|\, \int ^{\infty}_{-\infty} 
  \left( \rho_\delta \ast \mathrm{sgn} \right)(\phi) \, 
  \bigl(\, \partial _x U + \partial _x U^r \, \bigr)^{p-1} 
  \, \partial _x^2 U^r
  \, \mathrm{d}x \, \biggr|\, (t) \\
& \qquad 
  \leq C_{p} (1+t)^{-\frac{2p}{p+1}}
\quad \bigl( \, t \ge T_0 \, \bigr). 
\end{aligned}
\end{align}
Then, substituting (5.8), (5.9), (5.10) and (5.11) into (5.7), 
we have the desired $L^1$-estimate (5.1). 

\medskip
Next, we show 
the time-weighted $L^q$-energy estimates to $\phi$. 

\medskip

\noindent
{\bf Proposition 5.2.}\quad {\it
Suppose the same assumptions in Theorem 3.2. 
For any $q \in [ \, 1, \infty)$, 
there exist positive constants $\alpha$ and $C_{\alpha,p,q}$, 
such that the unique global solution in time $\phi$ of the Cauchy problem {\rm(3.6)} 
satisfies the following $L^q$-energy estimate 
 \begin{align}
 \begin{aligned}
 &(1+t)^\alpha  \| \,\phi(t) \,\|_{L^q}^q 
  + \int ^t_0 (1+\tau )^\alpha  G_{q}(\tau) 
    \, \mathrm{d}\tau \\
 & \quad 
  + \int ^t_0 (1+\tau )^\alpha  \int _{-\infty}^{\infty} 
    | \, \phi \, |^{q-2} 
    \bigl( \, \partial _x \phi \, \bigr)^2 \\
 & \qquad \qquad \qquad \qquad \quad \; \, 
    \times 
    \left( \, 
    \bigl| \partial _x \phi \bigr|^{p-1} 
    + \bigl| \partial _x U \bigr|^{p-1} 
    + \bigl| \partial _x U^r \bigr|^{p-1}  
    \, \right) \, \mathrm{d}x \mathrm{d}\tau  \\
 & \quad 
  + \int ^t_0 (1+\tau )^\alpha  \int _{-\infty}^{\infty} 
    | \, \phi \, |^{q-2} \, 
    \left| \, 
    \bigl| \partial _x \phi + \partial _x U + \partial _x U^r \bigr|^{p-1} 
    - \bigl| \partial _x U + \partial _x U^r \bigr|^{p-1}  
    \, \right| \\
 & \qquad \qquad \qquad \qquad \quad \; \, 
   \times 
    \left| \, 
    \bigl( \, \partial _x \phi + \partial _x U + \partial _x U^r \, \bigr)^2 
    - \bigl( \, \partial _x U + \partial _x U^r \, \bigr)^2 \, \right| 
    \, \mathrm{d}x \mathrm{d}\tau  \\
 &\leq C_{\alpha,p,q} \| \,\phi_0 \,\|_{L^q}^q 
       + C\left( \, \alpha, \, p, \, q, \, \phi_0 \, \right) \, 
       (1+t)^{\alpha - \frac{q-1}{2p}}
       \quad \bigl( \, t \ge T_0 \, \bigr). 
 \end{aligned}
 \end{align}
}

\medskip

The proof of Proposition 5.2 is given by the
following two lemmas.

\medskip

\noindent
{\bf Lemma 5.3.}\quad {\it
For any 
$1\leq q < \infty$, 
there exist positive constants $\alpha$ and $C_q$ such that 
\begin{align}
\begin{aligned}
&(1+t)^\alpha  \| \,\phi(t) \,\|_{L^q}^q 
  + q \, (q-1)
    \int ^t_0 (1+\tau )^\alpha  G_{q}(\tau) 
    \, \mathrm{d}\tau \\
 & 
  + C_{q} 
    \int ^t_0 (1+\tau )^\alpha 
    \int _{-\infty}^{\infty} 
    | \, \phi \, |^{q-2} 
    \bigl( \, \partial _x \phi \, \bigr)^2 \\
 & \qquad \qquad \qquad \qquad \quad \quad \quad \; \, 
   \times 
    \left( \, 
    \bigl| \partial _x \phi \bigr|^{p-1} 
    + \bigl| \partial _x U \bigr|^{p-1}  
    + \bigl| \partial _x U^r \bigr|^{p-1}  
    \, \right) \, \mathrm{d}x \mathrm{d}\tau  \\
 & 
  + C_{q} 
    \int ^t_0 (1+\tau )^\alpha  \int _{-\infty}^{\infty} 
    | \, \phi \, |^{q-2} \, 
    \left| \, 
    \bigl| \partial _x \phi + \partial _x \tilde{U} \bigr|^{p-1} 
    - \bigl| \partial _x \tilde{U} \bigr|^{p-1}  
    \, \right| \\
 & \qquad \qquad \qquad \qquad \quad \quad \quad \; \, 
   \times 
    \left| \, 
    \bigl( \, \partial _x \phi + \partial _x \tilde{U} \, \bigr)^2 
    - \bigl( \, \partial _x \tilde{U} \, \bigr)^2 \, \right| 
    \, \mathrm{d}x \mathrm{d}\tau  \\
&\leq \| \,\phi_0 \,\|_{L^q}^q 
      + \alpha 
      \int ^t_0 (1+\tau )^{\alpha -1} 
      \| \,\phi(\tau ) \,\|_{L^p}^p 
      \, \mathrm{d}\tau \\
&\qquad \quad \quad \; \, \,  
 + q \int ^t_0 (1+\tau )^\alpha 
   \| \,\phi(\tau) \,\|_{L^\infty}^{q-1}
   \left| \left| \, 
    \widetilde{F_p}(U,U^r) (\tau ) \, \right| \right|_{L^1} 
    \, \mathrm{d}\tau \\
&\qquad \quad \quad \; \, \, 
 + C_{q} \int ^t_0 (1+\tau )^\alpha 
   \widetilde{\widetilde{F_p}}(\phi, U, U^r) (\tau)
   \, \mathrm{d}\tau 
   \quad \bigl( \, t \ge T_0 \, \bigr), 
\end{aligned}
\end{align}
where 
\begin{align*}
\begin{aligned}
& \widetilde{\widetilde{F_p}}(\phi, U, U^r) (t) \\
&:= \left\{\begin{array} {ll} 
   \displaystyle{ 
   \| \,\phi(t) \,\|_{L^\infty}^{q-1} \, 
   \biggl|\, \int ^{\infty}_{-\infty} 
   \partial _x \Bigl( \, 
   \bigl(\, \partial _x U + \partial _x U^r \, \bigr)^{p} 
   - \bigl(\, \partial _x U \, \bigr)^{p} \, 
   \Bigr) \, \mathrm{d}x \, \biggr| \, (t)}
   \quad \bigl( \, 1<q<2 \, \bigr), \\[15pt]
   \displaystyle{ 
   \| \,\phi(t) \,\|_{L^\infty}^{q-2} \, 
   \biggl|\, \int ^{\infty}_{-\infty} 
   \Bigl( \, 
   \bigl(\, \partial _x U + \partial _x U^r \, \bigr)^{p} 
   - \bigl(\, \partial _x U \, \bigr)^{p} \, 
   \Bigr)^{\frac{p+1}{p}} \, \mathrm{d}x \, \biggr| \, (t)
   \quad \quad \; \: \, \bigl( \, q \geq 2 \, \bigr). 
   }
   \end{array}
   \right.\, 
\end{aligned}
\end{align*}
}

\bigskip

\noindent
{\bf Lemma 5.4.}\quad {\it
Assume $p>1$ and $1\leq q < \infty$. 
We have the following interpolation inequalities.  

 \noindent
 {\rm (1)}\ \ For any $1\leq r < \infty$, there exists a positive 
 constant $C_{p,q,r}$  such that 
 \begin{align*}
 \begin{aligned}
 &\Vert \, \phi (t) \, \Vert _{L^r }
  \leq 
  C_{p,q,r} \left( \, \int _{-\infty}^{\infty} | \, \phi \, |^2 \, \mathrm{d}x \, \right)
  ^{\frac{pr+p+q-1}{(2p+q-1)r}} \\
 & \qquad \quad \quad \quad \quad \; \: \, 
         \times 
         \left( \, \int _{-\infty}^{\infty} | \, \phi \, |^{q-2} 
         \bigl| \, \partial _x \phi \, \bigr|^{p+1} \, \mathrm{d}x \, \right)
         ^{\frac{r-1}{(2p+q-1)r}} \quad \bigl( t \ge 0 \bigr). 
 \end{aligned}
 \end{align*}

 \noindent
 {\rm (2)}\ \ 
 There exists a positive 
 constant $C_{p,q}$ 
 such that 
 \begin{align*}
 \begin{aligned}
 &\Vert \, \phi (t)  \, \Vert _{L^\infty }
  \leq 
  C_{p,q} \left( \, \int _{-\infty}^{\infty} | \, \phi \, |^2 \, \mathrm{d}x \, \right)
  ^{\frac{p}{2p+q-1}} \\
 & \qquad \quad \quad \quad \quad \; \, 
         \times 
         \left( \, \int _{-\infty}^{\infty} | \, \phi \, |^{q-2} 
         \bigl| \, \partial _x \phi \, \bigr|^{p+1} \, \mathrm{d}x \, \right)
         ^{\frac{1}{2p+q-1}} \quad \bigl( t \ge 0 \bigr). 
 \end{aligned}
 \end{align*}
}

\medskip

\noindent
The proofs of Lemma 5.3, Lemma 5.4 and Proposition 5.2 are given in 
the quite similar way as those of Lemma 4.1, Lemma 4.2 and Proposition 4.1, 
so we omit them. 
We particularly note that we have by Proposition 5.2 
\begin{equation}
\| \,\phi (t) \,\|_{L^q}
\leq C( \, p, \, q, \, \phi_0 \, )\, (1+t)^
     {-\frac{1}{2p}\left( 1- \frac{1}{q} \right)} 
\end{equation}
for $1\leq q < \infty$. 

\bigskip

We shall finally obtain the time-decay estimates 
for the higher order derivatives, that is, 
$\partial _x \phi$ and $\partial _x u$, 
and also get the $L^{\infty}$-estimate for $\phi$. 

\medskip

\noindent
{\bf Proposition 5.3.}\quad {\it
Suppose the same assumptions in Theorem 3.2. 
There exist positive constants $\alpha$ and $C_{\alpha,p}$, 
such that the unique global solution in time $\phi$ of the Cauchy problem {\rm(3.6)} 
satisfies the following $L^{p+1}$-energy estimate 
\begin{align}
\begin{aligned}
 &(1+t)^\alpha  
  \| \, \partial _x u (t)\, \| _{L^{p+1}}^{p+1} 
  + 
    \int ^t_0 (1+\tau )^\alpha  \int _{-\infty}^{\infty} 
    \bigl| \, \partial _x u \, \bigr|^{2(p-1)} 
    \left( \, \partial _x^2 u \, \right)^2 
    \, \mathrm{d}x \mathrm{d}\tau \\
 & \qquad \qquad \qquad \quad \quad \quad \, 
  + \int ^t_0 (1+\tau )^\alpha  
    \| \, \partial _x u (\tau) \, \| _{L^{p+2}}^{p+2} 
    \, \mathrm{d}\tau  \\
 &\leq C_{\alpha,p} 
       \| \, \partial _x u_{0} \, \| _{L^{p+1}}^{p+1} 
       + C\left( \, \alpha, \, p, \, \phi_0, \, \partial _x u_{0} \, \right) \, 
       (1+t)^{\alpha - \frac{1}{2p}}
       \quad \bigl( \, t \ge T_0 \, \bigr). 
\end{aligned}
\end{align}
}

\medskip

\noindent
To obtain Proposition 5.3, we first show the following. 

\medskip

\noindent
{\bf Lemma 5.5.}\quad {\it
It follows that 
\begin{align}
\begin{aligned}
 &(1+t)^\alpha  
  \| \, \partial _x u (t)\, \| _{L^{p+1}}^{p+1} \\
 & \quad 
  + \mu \, p^2 \, (p+1) 
    \int ^t_0 (1+\tau )^\alpha  \int _{-\infty}^{\infty} 
    \bigl| \, \partial _x u \, \bigr|^{2(p-1)} 
    \left( \, \partial _x^2 u \, \right)^2 
    \, \mathrm{d}x \mathrm{d}\tau \\
 & \qquad \quad \quad \; \; \; \: \, 
  + p \int ^t_0 (1+\tau )^\alpha  
    \int _{\partial _x u \geq 0} 
    f''(u) \left|\, \partial _x u \, \right|^{p+2} 
    \, \mathrm{d}x \mathrm{d}\tau  \\
 &= \| \, \partial _x u_{0} \, \| _{L^{p+1}}^{p+1} 
    + \alpha \int ^t_0 (1+\tau )^{\alpha -1} 
      \| \, \partial _x u (\tau) \, \| _{L^{p+1}}^{p+1} 
      \, \mathrm{d}\tau \\
 & \qquad \quad \quad \; \; \; \: \, 
  + p \int ^t_0 (1+\tau )^\alpha 
    \int _{\partial _x u < 0} 
    f''(u) \left|\, \partial _x u \, \right|^{p+2} 
    \, \mathrm{d}x \mathrm{d}\tau  
       \quad \bigl( \, t \ge T_0 \, \bigr). 
\end{aligned}
\end{align}
}

\medskip

\noindent
{\bf Proof of Lemma 5.5.}\quad
Multiplying the equation in the problem (1.1), that is, 
$$
\partial_tu +\partial_x \bigl(f(u) \bigr)
  = \mu \, 
    \partial_x \left( \, 
    \left| \, \partial_xu \, \right|^{p-1} \partial_xu \, 
    \right)
$$
by 
$$
- \partial_x \left( \, 
\left| \, \partial_xu \, \right|^{p-1} \partial_xu \, 
\right), 
$$ 
we obtain the divergence form 
\begin{align}
\begin{aligned}
&\partial_t 
 \left(\frac{1}{p+1} \left|\, \partial _x u \, \right|^{p+1} \right) 
 + \partial _x \Bigl( \, 
   - \, \left|\, \partial _x u \, \right|^{p-1} 
   \partial _x u \cdot \partial _t u \, \Bigr) \\
&\qquad \; \, 
 + \partial _x \left( 
 - \, \frac{p}{p+1} \, 
 f'(u) \left|\, \partial _x u \, \right|^{p+1} 
 \, \right) \\
&\qquad \; \, + \frac{p}{p+1} \, 
          f''(u) \left|\, \partial _x u \, \right|^{p+1}\partial _x u 
+\mu \, p\, q \left|\, \partial _x u \, \right|^{2(p-1)} 
   \bigl( \, \partial _x^2 u \, \bigr)^2 
= 0. 
\end{aligned}
\end{align}
Integrating the divergence form (5.17) with respect to $x$, 
we have 
\begin{align}
\begin{aligned}
&\frac{1}{p+1} \, 
\frac{\mathrm{d}}{\mathrm{d}t}\, 
\,\Vert \, \partial _x u (t) \, \Vert_{L^{p+1}}^{p+1} 
+\mu \, p^2 \int ^{\infty }_{-\infty } 
\left|\, \partial _x u \, \right|^{2(p-1)} 
\bigl( \, \partial _x^2 u \, \bigr)^2 
\, \mathrm{d}x \\
&\qquad \qquad \qquad \qquad \quad \; \; \; \: 
 +\frac{p}{p+1} \int ^{\infty }_{-\infty } 
  f''(u) \left|\, \partial _x u \, \right|^{p+1}\partial _x u 
  \, \mathrm{d}x 
= 0. 
\end{aligned}
\end{align}
We separate the integral region 
to the third term on the left-hand side of (5.18) as 
\begin{align}
\begin{aligned}
&\int ^{\infty }_{-\infty } 
 f''(u) \left|\, \partial _x u \, \right|^{p+1}\partial _x u 
 \, \mathrm{d}x \\
&= \int _{\partial _x u \geq 0 } + \int _{\partial _x u < 0 } \\
&= \int _{\partial _x u \geq 0 } 
   f''(u) \left|\, \partial _x u \, \right|^{p+2} \, \mathrm{d}x 
  - \int _{\partial _x u < 0 } 
  f''(u) \left|\, \partial _x u \, \right|^{p+2} \, \mathrm{d}x. 
\end{aligned}
\end{align}
Substituting (5.19) into (5.18), we get the following equality 
\begin{align}
\begin{aligned}
&\frac{1}{p+1} \, 
\frac{\mathrm{d}}{\mathrm{d}t}\, 
\,\Vert \, \partial _x u (t) \, \Vert_{L^{p+1}}^{p+1} 
+\mu \, p^2 \int ^{\infty }_{-\infty } 
\left|\, \partial _x u \, \right|^{2(p-1)} 
\bigl( \, \partial _x^2 u \, \bigr)^2 
\, \mathrm{d}x \\
& \; \: \, +\frac{p}{p+1} \int _{\partial _x u \geq 0 } 
   f''(u) \left|\, \partial _x u \, \right|^{p+2} \, \mathrm{d}x 
= \frac{p}{p+1} \int _{\partial _x u < 0 } 
  f''(u) \left|\, \partial _x u \, \right|^{p+2} \, \mathrm{d}x. 
\end{aligned}
\end{align}
Multiplying (5.20) by 
$(1+t)^{\alpha}$ with $\alpha>0$ 
and integrating over $(0,t)$ with respect to the time, 
we complete the proof of Lemma 5.5. 

\bigskip

\noindent
{\bf Proof of Proposition 5.3.}\quad
We use the following important results (cf. \cite{yoshida'}). 

\medskip

\noindent
{\bf Lemma 5.6.}\quad {\it
For any $s \geq 0$, there exists a positive constant $C_{s}$ 
such that 
\begin{align}
\begin{aligned}
\int _{\partial _x u < 0 } 
f''(u) \left|\, \partial _x u \, \right|^{s} \, \mathrm{d}x 
\leq C_{s}
     \int _{\partial _x u < 0 } 
     \left|\, \partial _x \phi \, \right|^{s} \, \mathrm{d}x. 
\end{aligned}
\end{align}
}

\noindent
In fact, taking care of the relation by using Lemma 2.2 and Lemma 2.3
\begin{align}
\begin{aligned}
\partial _x u = \partial _x \tilde{U} + \partial _x \phi <0 \, 
\Longleftrightarrow  \, \partial _x \phi <0, 
                     \, \partial _x \tilde{U} < \bigl|\, \partial _x \phi \, \bigr|, 
\end{aligned}
\end{align}
we immediately have 
\begin{align}
\begin{aligned}
\int _{\partial _x u < 0 } 
&f''(u) \left|\, \partial _x u \, \right|^{s} \, \mathrm{d}x \\
&\leq 2^{s} \, 
     \left( \, \max_{ | u | \leq \widetilde{C}} f''(u) \, \right)
     \int _{\partial _x \phi < 0, \partial _x \tilde{U} < | \partial _x \phi | } 
     \left|\, \partial _x \phi \, \right|^{s} \, \mathrm{d}x.
\end{aligned}
\end{align}

\medskip
Since $ \partial _x u $ is absolutely continuous, 
we first note that 
for any 
$
x \in \bigl\{ \, x \in \mathbb{R} \, \, \bigr. 
 \bigl| \, \partial _x u \, < \, 0 \, \bigr\}, 
$
there exsists 
$
x_{k} \in \mathbb{R}\cup \{ - \infty \}
$
such that 
$$
\partial _x u (x_{k}) = 0, \; \; 
  \partial _x u (y) \, < \, 0 \; 
  \bigl( \, y \in 
  ( x_{k},x ) \, \bigr).
$$
Therefore, 
it follows that for such $x$ and $x_{k}$ with $q\geq p\, (\, >1\, )$, 
\begin{align}
\left|\, \partial _x u \, \right|^q 
= \left(\, - \partial _x u \, \right)^q 
= q \int_{x_k}^{x} 
   \left(\, - \partial _x u \, \right)^{q-1} \left(\, - \partial _x^2 u \, \right) 
   \, \mathrm{d}y 
\end{align}
By using the Cauchy-Schwarz inequality, we have 

\medskip

\noindent
{\bf Lemma 5.7.}\quad {\it
It holds that 
\begin{align}
\begin{aligned}
& \int _{\partial _x u < 0 } 
  \left|\, \partial _x u \, \right|^{p+2} \, \mathrm{d}x \\
& \le C_{p} 
     \left( \, 
     \int _{\partial _x u < 0 } 
     \left|\, \partial _x u \, \right|^{2(p-1)} 
     \bigl( \, \partial _x^2 u \, \bigr)^2 
     \, \mathrm{d}x \, 
     \right)^{\frac{1}{3p+1}} 
     \left( \, 
     \int _{\partial _x u < 0 } 
     \left|\, \partial _x u \, \right|^{p+1} 
     \, \mathrm{d}x \, 
     \right)^{\frac{3p+2}{3p+1}}. 
\end{aligned}
\end{align}
}
By using Young's inequality to (5.25), we also have 

\medskip

\noindent
{\bf Lemma 5.8.}\quad {\it
It follows that for any $\epsilon>0$, 
there exists a positive constant $C_p({\epsilon})$ such that, 
\begin{align}
\begin{aligned}
& \int _{\partial _x u < 0 } 
  \left|\, \partial _x u \, \right|^{p+2} \, \mathrm{d}x \\
& \le \epsilon 
      \int _{\partial _x u < 0 } 
      \left|\, \partial _x u \, \right|^{2(p-1)} 
      \bigl( \, \partial _x^2 u \, \bigr)^2 
      \, \mathrm{d}x \, 
      + C_p({\epsilon}) \left( \, 
        \int _{\partial _x u < 0 } 
        \left|\, \partial _x u \, \right|^{p+1} 
        \, \mathrm{d}x \, 
        \right)^{\frac{3p+2}{3p}}. 
\end{aligned}
\end{align}
}
By using Lemma 5.6, Lemma 5.7 and Lemma 5.8 
with $\epsilon = \frac{\mu \, p^2 \, (p+1)}{2}$, we have 
\begin{align}
\begin{aligned}
 &(1+t)^\alpha  
  \| \, \partial _x u (t)\, \| _{L^{p+1}}^{p+1} \\
 & \quad 
  + \frac{\mu \, p^2 \, (p+1)}{2} 
    \int ^t_0 (1+\tau )^\alpha  \int _{-\infty}^{\infty} 
    \bigl| \, \partial _x u \, \bigr|^{2(p-1)} 
    \left( \, \partial _x^2 u \, \right)^2 
    \, \mathrm{d}x \mathrm{d}\tau \\
 & \qquad \quad \quad \; \; \; \: \, 
  + p \int ^t_0 (1+\tau )^\alpha  
    \int _{\partial _x u \geq 0} 
    f''(u) \left|\, \partial _x u \, \right|^{p+1} 
    \, \mathrm{d}x \mathrm{d}\tau  \\
 &\leq \| \, \partial _x u_{0} \, \| _{L^{p+1}}^{p+1} \\
 & \quad 
       + \alpha \int ^t_0 (1+\tau )^{\alpha -1} 
         \Bigl( \, 
         \| \, \partial _x \phi (\tau) \, \| _{L^{p+1}}^{p+1} 
         + \| \, \partial _x \tilde{U}(\tau) \, \| _{L^{p+1}}^{p+1} 
         \, \Bigr) 
         \, \mathrm{d}\tau \\ 
 & \quad 
  + C_{p} \int ^t_0 (1+\tau )^\alpha 
    \left( \, 
  \int _{\partial _x u < 0 } 
  \left|\, \partial _x u \, \right|^{p+1} 
  \, \mathrm{d}x \, 
  \right)^{\frac{2}{3p}+1} \mathrm{d}\tau . 
\end{aligned}
\end{align}
By using Proposition 5.2, we get the following time-decay estimates. 

\medskip

\noindent
{\bf Lemma 5.9.}\quad {\it
There exist positive constants $\alpha \gg 1$ and $C_{\alpha,p,q}$, 
such that 
 \begin{align}
 \begin{aligned}
 &\int ^t_0 (1+\tau )^\alpha  \int _{-\infty}^{\infty} 
    | \, \phi \, |^{q-2} 
    \bigl( \, \partial _x \phi \, \bigr)^2 
    \left( \, 
    \bigl| \partial _x \phi \bigr|^{p-1} 
    + \bigl| \partial _x U \bigr|^{p-1}  
    + \bigl| \partial _x U^r \bigr|^{p-1}  
    \, \right) \, \mathrm{d}x \mathrm{d}\tau  \\
 &\leq C\left( \, \alpha, \, p, \, q, \, \phi_0 \, \right) \, 
       (1+t)^{\alpha - \frac{q-1}{2p}}
       \quad \bigl( \, t \ge T_0 \, \bigr).
 \end{aligned}
 \end{align}
}

\medskip

\noindent
By using Lemma 5.9 
with $\alpha \mapsto \alpha - 1 \gg 1$ and $q=2$, we have 
\begin{align}
\alpha \int ^t_0 (1+\tau )^{\alpha - 1} 
\| \, \partial _x \phi (\tau) \, \| _{L^{p+1}}^{p+1} 
\, \mathrm{d}\tau  
\leq C\left( \, \alpha, \, p, \, \phi_0 \, \right) \, 
     (1+t)^{\alpha - \frac{2p+1}{2p}}. 
\end{align}
We can also estimate by using Lemma 2.2 and Lemma 2.3 as 
\begin{align}
\alpha \int ^t_0 (1+\tau )^{\alpha - 1} 
\| \, \partial _x \tilde{U} (\tau) \, \| _{L^{p+1}}^{p+1} 
\, \mathrm{d}\tau  
\leq C\left( \, \alpha, \, p \, \right) \, 
     (1+t)^{\alpha - \frac{p}{p+1}}. 
\end{align}
By using the uniform boundedness in Lemma 3.2, that is, 
$$
\| \, \partial _x u (t) 
\, \| _{L^{p+1}}^{p+1} 
\leq 
 C_{p}\bigl(\, \| \, \phi _0 \, \| _{L^2}, 
 \| \, \partial _x u_0 \, \| _{L^{p+1}} \, \bigr) 
$$
and Lemma 5.9 
with $q=2$, we can estimate as 
\begin{align}
\begin{aligned}
&C_{p} \int ^t_0 (1+\tau )^\alpha 
 \left( \, 
 \int _{\partial _x u < 0 } 
 \left|\, \partial _x u \, \right|^{p+1} 
 \, \mathrm{d}x \, 
 \right)^{\frac{2}{3p}+1} \mathrm{d}\tau \\
&\leq C_{p} \int ^t_0 (1+\tau )^\alpha 
      \int _{\partial _x u < 0 } 
      \left|\, \partial _x \phi \, \right|^{p+1} 
      \, \mathrm{d}x 
      \cdot 
      \| \, \partial _x u (\tau) \, \| _{L^{p+1}}^{\frac{2(p+1)}{3p}} 
      \, \mathrm{d}\tau \\
&\leq C\left( \, p, \, \phi_0, \, \partial _x u_0 \, \right) \, 
      \int ^t_0 (1+\tau )^\alpha 
      \int _{-\infty}^{\infty} 
      \left|\, \partial _x \phi \, \right|^{p+1} 
      \, \mathrm{d}x \mathrm{d}\tau \\
&\leq C\left( \, \alpha, \, p, \, \phi_0, \, \partial _x u_0 \, \right) \, 
      (1+t)^{\alpha - \frac{1}{2p}}. 
\end{aligned}
\end{align}
Substituting (5.29), (5.30) and (5.31) into (5.27), 
we complete the proof of Proposition 5.3. 
In particular, we have 
\begin{align}
\| \, \partial _x u (t) \, \| _{L^{p+1}}^{p+1} 
\, \mathrm{d}\tau  
\leq C\left( \, p, \, \phi_0, \, \partial _x u_0 \, \right) \, 
     (1+t)^{- \frac{1}{2p}}, 
\end{align}
and 
\begin{align}
\begin{aligned}
\| \, \partial _x \phi (t) \, \| _{L^{p+1}}^{p+1} 
\, \mathrm{d}\tau  
& \leq \| \, \partial _x u (t) \, \| _{L^{p+1}}^{p+1} 
       + \| \, \partial _x U^r (t) \, \| _{L^{p+1}}^{p+1} \\
& \leq C\left( \, p, \, \phi_0, \, \partial _x u_0 \, \right) \, 
       (1+t)^{- \frac{1}{2p}} 
\end{aligned}
\end{align}
for $1 \leq q < \infty$.

\medskip

{\bf Proof of Theorem 3.2.}\ 
We already have proved the decay estimate 
of $\| \,\phi(t) \,\|_{L^q}$ with $1 \le q < \infty$. 
Therefore we only 
show 
the following time-decay estimate for the higher order derivative
\begin{align}
\begin{aligned}
&\bigl|\bigl|\,
 \partial _x u(t) \,
 \bigr|\bigr|_{L^{p+1} }, \, \; \; 
 \left|\left|\,
 \partial _x \phi (t) 
 \,\right|\right|_{L^{p+1} } \\[5pt]
& \leq \left\{\begin{array} {ll} 
      C( \, \epsilon, \, p, \, \phi_0, \, \partial _x u_0 \, ) \, 
      (1+t)^{-\frac{p}{(p+1)^2}}\\[15pt] 
      \,  \, \, \: \; \; \quad \qquad 
      \left( \, 1 < p \le 
      \displaystyle{\frac{7}{12} + 
               \sqrt{\frac{73}{144} - \frac{(p+1)^2(3p-2)}{3}\, \epsilon } }
      \, \right),\\[15pt] 
      C( \, \epsilon, \, p, \, \phi_0, \, \partial _x u_0 \, ) \, 
      (1+t)^{-\frac{3}{2(p+1)(3p-2)} + \epsilon}\\[15pt] 
      \,  \, \, \; \; \quad \qquad \qquad 
      \left( \,  
      \displaystyle{\frac{7}{12} + 
               \sqrt{\frac{73}{144} - \frac{(p+1)^2(3p-2)}{3}\, \epsilon } }
      < p \, \right)
      \end{array}
      \right.\,
\end{aligned}
\end{align}
for any $0<\epsilon \ll 1$, 
and the $L^{\infty}$-estimate for $\phi$. 

\medskip

We first prove (5.34). 
Substituting (5.32) into (5.31), 
we have 
\begin{align}
\begin{aligned}
&C_{p} \int ^t_0 (1+\tau )^\alpha 
      \int _{\partial _x u < 0 } 
      \left|\, \partial _x \phi \, \right|^{p+1} 
      \, \mathrm{d}x 
      \cdot 
      \| \, \partial _x u (\tau) \, \| _{L^{p+1}}^{\frac{2(p+1)}{3p}} 
      \, \mathrm{d}\tau \\
&\leq C\left( \, p, \, \phi_0, \, \partial _x u_0 \, \right) \, 
      \int ^t_0 (1+\tau )
      ^{\alpha - \frac{1}{2p} \cdot \frac{2}{3p}}
      \int _{-\infty}^{\infty} 
      \left|\, \partial _x \phi \, \right|^{p+1} 
      \, \mathrm{d}x \mathrm{d}\tau. 
\end{aligned}
\end{align}

\medskip

\noindent
By using Lemma 5.9 
with $\alpha \mapsto \alpha - \frac{1}{2p} \cdot \frac{2}{3p} \gg 1$ 
and $q=2$, we also have 
\begin{align}
\alpha \int ^t_0 (1+\tau )^{\alpha - 1} 
\| \, \partial _x \phi (\tau) \, \| _{L^{p+1}}^{p+1} 
\, \mathrm{d}\tau  
\leq C\left( \, \alpha, \, p, \, \phi_0 \, \right) \, 
     (1+t)^{\alpha - \frac{1}{2p} \cdot \frac{2}{3p} - \frac{1}{2p}}. 
\end{align}
Substituting (5.36) into (5.27), 
we have 
\begin{align}
\begin{aligned}
&\bigl|\bigl|\,
 \partial _x u(t) \,
 \bigr|\bigr|_{L^{p+1} }^{p+1}, \, \; \; 
 \left|\left|\,
 \partial _x \phi (t) 
 \,\right|\right|_{L^{p+1} }^{p+1} \\[5pt]
& \leq C( \, p, \, \phi_0, \, \partial _x u_0 \, ) \\
& \quad \times 
       \left( \, 
       (1+t)^{-\frac{2p+1}{2p}} + (1+t)^{-\frac{p}{p+1}}  
       + (1+t)
       ^{-\left( \frac{1}{2p} \cdot \frac{2}{3p} + \frac{1}{2p} \right)} 
       \,  \right) \\[5pt] 
& \leq C( \, p, \, \phi_0, \, \partial _x u_0 \, ) \, 
       \left( \, 
       (1+t)^{-\frac{p}{p+1}} 
       + (1+t)
       ^{-\left( \frac{1}{2p} \cdot \frac{2}{3p} + \frac{1}{2p} \right)} 
       \,  \right) 
\end{aligned}
\end{align}
Iterating ``$\infty$''-times the above process, we will get 
\begin{align}
\begin{aligned}
&\bigl|\bigl|\,
 \partial _x u(t) \,
 \bigr|\bigr|_{L^{p+1} }^{p+1}, \, \; \; 
 \left|\left|\,
 \partial _x \phi (t) 
 \,\right|\right|_{L^{p+1} }^{p+1} \\[5pt]
& \leq C( \, \epsilon, \, p, \, \phi_0, \, \partial _x u_0 \, ) \, 
       \left( \, 
       (1+t)^{-\frac{p}{p+1}} 
       + (1+t)
         ^{- \mathlarger{ \frac{1}{2p}} 
         \substack{{\infty }\\{\substack{{\lsum }\\{n=0}}}}
         \mathlarger{\left( \frac{2}{3p} \right)^{n}} + \mathlarger{\epsilon } }
       \,  \right) \\[15pt] 
& \leq C( \, \epsilon, \, p, \, \phi_0, \, \partial _x u_0 \, ) \, 
       \left( \, 
       (1+t)^{-\frac{p}{p+1}} 
       + (1+t)^{- \frac{3}{2(3p-2)} + \epsilon }
       \,  \right) \\[15pt] 
& \leq \left\{\begin{array} {ll} 
      C( \, \epsilon, \, p, \, \phi_0, \, \partial _x u_0 \, ) \, 
      (1+t)^{-\frac{p}{p+1}}\\[15pt] 
      \, \: \; \; \qquad \qquad 
      \left( \, 1 < p 
      \le \displaystyle{\frac{7}{12} + 
               \sqrt{\frac{73}{144} - \frac{(p+1)(3p-2)}{3}\, \epsilon } } \, \right),\\[15pt] 
      C( \, \epsilon, \, p, \, \phi_0, \, \partial _x u_0 \, ) \, 
      (1+t)^{-\frac{3}{2(3p-2)} + \epsilon}\\[15pt]
      \, \: \; \quad \quad \qquad \qquad 
      \left( \,  
      \displaystyle{\frac{7}{12} + 
               \sqrt{\frac{73}{144} - \frac{(p+1)(3p-2)}{3}\, \epsilon } }
      < p \, \right)
      \end{array}
      \right.\,
\end{aligned}
\end{align}
for any $0<\epsilon \ll 1$. 

\medskip

\noindent
Thus, we get (5.34). 

We finally show the $L^{\infty}$-estimate for $\phi$ 
by using the Gagliardo-Nirenberg inequality.  
Substituting (5.14) and (5.34) into (4.17), we get 
\begin{align}
\begin{aligned}
&\bigl|\bigl|\,
 \phi(t) \,
 \bigr|\bigr|_{L^{\infty } } \\[5pt] 
&\leq \left\{\begin{array} {ll} 
      C( \, \epsilon, \, p, \, \theta, \, \phi_0, \, \partial _x u_0 \, ) \, 
      (1+t)
      ^{-\frac{1}{2p} + 
      \left( \frac{2p+1}{2p(p+1)} - \frac{p}{(p+1)^2} \right) \theta} \\[15pt] 
      \, \: \: \; \; \quad \qquad \qquad \qquad \, 
      \left( \, 1 < p \le 
      \displaystyle{\frac{7}{12} + 
               \sqrt{\frac{73}{144} - \frac{(p+1)^2(3p-2)}{3}\, \epsilon } }
      \, \right),\\[10pt] 
      C( \, \epsilon, \, p, \, \theta, \, \phi_0, \, \partial _x u_0 \, ) \, 
      (1+t)
      ^{-\frac{1}{2p} + 
      \left( \frac{2p+1}{2p(p+1)} - \frac{3}{2(p+1)(3p-2)}  + \epsilon \right) \theta}\\[15pt] 
      \, \: \: \; \; \quad \qquad \qquad \qquad \qquad 
      \left( \,  
      \displaystyle{\frac{7}{12} + 
               \sqrt{\frac{73}{144} - \frac{(p+1)^2(3p-2)}{3}\, \epsilon } }
      < p \, \right)
      \end{array}
      \right.\, 
\end{aligned}
\end{align}
for $\theta \in (0,1\, ]$ and any $0<\epsilon \ll 1$. 
Consequently, we do complete the proof of Theorem 3.2.
\bigskip 

\noindent
\section{$L^{r+1}$-estimate for the higher order derivative with $r>p$}
In this section, 
we show the time-decay estimates for the higher order derivative 
in the $L^{r+1}$-norm 
with $r>p$, 
in the case where $\phi_0 \in L^1 \cap L^2$ 
with $\partial _x u_0 \in L^{p+1} \cap L^{r+1}$, 
that is, Theorem 3.3. 

\medskip

\noindent
{\bf Proposition 6.1.}\quad {\it
Suppose the same assumptions in Theorem 3.3. 
For any $r>p$, 
there exist positive constants $\alpha$ and $C_{\alpha,p,r}$, 
such that the unique global solution in time $\phi$ of the Cauchy problem {\rm(3.6)} 
satisfies the following $L^{r+1}$-energy estimate 
\begin{align}
\begin{aligned}
 &(1+t)^\alpha  
  \| \, \partial _x u (t)\, \| _{L^{r+1}}^{r+1} 
  + \int ^t_0 (1+\tau )^\alpha  \int _{-\infty}^{\infty} 
    \bigl| \, \partial _x u \, \bigr|^{p+r-2} 
    \left( \, \partial _x^2 u \, \right)^2 
    \, \mathrm{d}x \mathrm{d}\tau \\
 & \qquad \qquad \qquad \quad \quad \quad \, 
  + \int ^t_0 (1+\tau )^\alpha  
    \| \, \partial _x u (\tau) \, \| _{L^{r+2}}^{r+2} 
    \, \mathrm{d}\tau  \\
 &\leq C_{\alpha,p,r} 
       \| \, \partial _x u_{0} \, \| _{L^{r+1}}^{r+1} \\[5pt] 
 & \; \; 
     + \left\{\begin{array} {ll} 
       C( \, \alpha, \, \epsilon, \, p, \, r, \, \phi_0, \, \partial _x u_0 \, ) \, 
       (1+t)^{\alpha-\frac{4p(r-p)+7p+3}{6p(p+1)}} \\[15pt] 
       \left( \, 
       1 < p \le 
       \displaystyle{\frac{7}{12} + 
               \sqrt{\frac{73}{144} - \frac{(p+1)(3p-2)}{3}\, \epsilon } }, \; 
       r > \displaystyle{\frac{-4p^2+7p+3}{2p}} > p 
       \, \right),\\[15pt] 
       C( \, \alpha, \, \epsilon, \, p, \, r, \, \phi_0, \, \partial _x u_0 \, ) \, 
       (1+t)^{\alpha-\frac{r}{p+1} 
       + \frac{2(r-p+1)}{3p} \epsilon} \\[15pt] 
       \left( \, 
       1 < p \le 
       \displaystyle{\frac{7}{12} + 
               \sqrt{\frac{73}{144} - \frac{(p+1)(3p-2)}{3}\, \epsilon } }, \; 
       p < r \leq \displaystyle{\frac{-4p^2+7p+3}{2p}} 
       \, \right),\\[15pt] 
       C( \, \alpha, \, \epsilon, \, p, \, r, \, \phi_0, \, \partial _x u_0 \, ) \, 
       (1+t)^{\alpha-\frac{p+2r}{2p(3p-2)}
       + \frac{2(r-p+1)}{3p} \epsilon} \\[15pt] 
       \, \, \: \; \; \quad \quad \qquad \qquad \qquad \qquad 
       \left( \, 
       \displaystyle{\frac{7}{12} + 
               \sqrt{\frac{73}{144} - \frac{(p+1)(3p-2)}{3}\, \epsilon } }
       <p 
       \, \right)
       \end{array}
       \right.\,
\end{aligned}
\end{align}
for $t \ge T_0$ and any $0<\epsilon \ll 1$. 
}

\medskip

The proof of Proposition 6.1 is given by the following 
three lemmas. 
Because the proofs of them are similar to those of 
Lemma 5.5, Lemma 5.6, Lemma 5.7 and Lemma 5.8, 
we state only here. 

\medskip

\noindent
{\bf Lemma 6.1.}\quad {\it
There exist positive constants 
$C_{p,r}$ and $C_{\alpha,p,r}$ such that 
\begin{align}
\begin{aligned}
 &(1+t)^\alpha  
  \| \, \partial _x u (t)\, \| _{L^{p+1}}^{p+1} \\
 & \quad 
  + \mu \, p\, r \, (r+1) 
    \int ^t_0 (1+\tau )^\alpha  \int _{-\infty}^{\infty} 
    \bigl| \, \partial _x u \, \bigr|^{p+r-2} 
    \left( \, \partial _x^2 u \, \right)^2 
    \, \mathrm{d}x \mathrm{d}\tau \\
 & \qquad \quad \quad \; \; \; \: \, 
  + r \int ^t_0 (1+\tau )^\alpha  
    \int _{\partial _x u \geq 0} 
    f''(u) \left|\, \partial _x u \, \right|^{r+2} 
    \, \mathrm{d}x \mathrm{d}\tau  \\
 &\leq \| \, \partial _x u_{0} \, \| _{L^{p+1}}\\
 & \quad \; \; \; \: \, 
       + C_{\alpha,p,r} \int ^t_0 
         (1+\tau )
         ^{\alpha -{\frac{2p+r+1}{3p+1}}} 
         \left( \, 
         \int _{-\infty}^{\infty} 
         \left|\, \partial _x u \, \right|^{p+1} 
         \, \mathrm{d}x \, 
         \right)^{\frac{p+2r+1}{3p+1}}
         \, \mathrm{d}\tau \\
 & \quad \; \; \; \: \, 
  + C_{p,r} \int ^t_0 (1+\tau )^\alpha 
    \left( \, 
         \int _{\partial _x u < 0 } 
         \left|\, \partial _x u \, \right|^{p+1} 
         \, \mathrm{d}x \, 
         \right)^{\frac{p+2r+2}{3p}}
         \, \mathrm{d}\tau  
       \quad \bigl( \, t \ge T_0 \, \bigr). 
\end{aligned}
\end{align}
}

\bigskip

\noindent
{\bf Lemma 6.2.}\quad {\it
Assume $p>1$ and $r>p$. 
We have the following interpolation inequalities.  

 \noindent
 {\rm (1)}\ \ There exists a positive 
 constant $C_{p,r}$  such that 
 \begin{align*}
 \begin{aligned}
 &\Vert \, \partial _x u (t) \, \Vert _{L^{r+1} }
 \leq 
  C_{p,r} \left( \, 
     \int _{-\infty}^{\infty}
     \left|\, \partial _x u \, \right|^{2(p-1)} 
     \bigl( \, \partial _x^2 u \, \bigr)^2 
     \, \mathrm{d}x \, 
     \right)^{\frac{r-p}{(2p+r+1)(r+1)}} \\
 & \: \; \; \quad \quad \quad \qquad \qquad \times 
     \left( \, 
     \int _{\partial _x u < 0 } 
     \left|\, \partial _x u \, \right|^{p+1} 
     \, \mathrm{d}x \, 
     \right)^{\frac{p+2r+1}{(2p+r+1)(r+1)}}
 \quad \bigl( t \ge 0 \bigr). 
 \end{aligned}
 \end{align*}

 \noindent
 {\rm (2)}\ \ 
 There exists a positive 
 constant $C_{p,r}$ 
 such that 
 \begin{align*}
 \begin{aligned}
 &\Vert \, \partial _x u (t) \, \Vert _{L^\infty }
 \leq 
  C_{p,r} \left( \, 
     \int _{-\infty}^{\infty} 
     \left|\, \partial _x u \, \right|^{2(p-1)} 
     \bigl( \, \partial _x^2 u \, \bigr)^2 
     \, \mathrm{d}x \, 
     \right)^{\frac{1}{2p+r+1}} \\
 & \: \quad \quad \quad \qquad \qquad \times 
     \left( \, 
     \int _{\partial _x u < 0 } 
     \left|\, \partial _x u \, \right|^{p+1} 
     \, \mathrm{d}x \, 
     \right)^{\frac{1}{2p+r+1}}
 \quad \bigl( t \ge 0 \bigr). 
 \end{aligned}
 \end{align*}
}

\bigskip

\noindent
{\bf Lemma 6.3.}\quad {\it
Assume $p>1$ and $r>p$. 
We have the following interpolation inequalities.  

 \noindent
 {\rm (1)}\ \ There exists a positive 
 constant $C_{p,r}$  such that 
 \begin{align*}
 \begin{aligned}
 &\Vert \, \partial _x u (t) \, \Vert 
  _{L^{r+1}_{x} \left( \{ \partial _x u < 0 \} \right) }
 \leq 
  C_{p,r} \left( \, 
     \int _{\partial _x u < 0 } 
     \left|\, \partial _x u \, \right|^{2(p-1)} 
     \bigl( \, \partial _x^2 u \, \bigr)^2 
     \, \mathrm{d}x \, 
     \right)^{\frac{r-p}{(2p+r+1)(r+1)}} \\
 & \, \: \: \quad \qquad \qquad \qquad \qquad \qquad \times 
     \left( \, 
     \int _{\partial _x u < 0 } 
     \left|\, \partial _x u \, \right|^{p+1} 
     \, \mathrm{d}x \, 
     \right)^{\frac{p+2r+1}{(2p+r+1)(r+1)}}
 \quad \bigl( t \ge 0 \bigr). 
 \end{aligned}
 \end{align*}

 \noindent
 {\rm (2)}\ \ 
 There exists a positive 
 constant $C_{p,r}$ 
 such that 
 \begin{align*}
 \begin{aligned}
 &\Vert \, \partial _x u (t) \, \Vert 
  _{L^{\infty}_{x} \left( \{ \partial _x u < 0 \} \right) }
 \leq 
  C_{p,r} \left( \, 
     \int _{\partial _x u < 0 } 
     \left|\, \partial _x u \, \right|^{2(p-1)} 
     \bigl( \, \partial _x^2 u \, \bigr)^2 
     \, \mathrm{d}x \, 
     \right)^{\frac{1}{2p+r+1}} \\
 & \, \, \: \: \; \qquad \qquad \qquad \qquad \qquad \times 
     \left( \, 
     \int _{\partial _x u < 0 } 
     \left|\, \partial _x u \, \right|^{p+1} 
     \, \mathrm{d}x \, 
     \right)^{\frac{1}{2p+r+1}}
 \quad \bigl( t \ge 0 \bigr). 
 \end{aligned}
 \end{align*}
}

\medskip

{\bf Proof of Proposition 6.1.}\ 
By using (5.34), 
we estimate the each terms on the right-hand side of (6.2) as 
\begin{align}
\begin{aligned}
&C_{\alpha,p,r} \int ^t_0 
         (1+\tau )
         ^{\alpha -{\frac{2p+r+1}{3p+1}}} 
         \left( \, 
         \int _{-\infty}^{\infty} 
         \left|\, \partial _x u \, \right|^{p+1} 
         \, \mathrm{d}x \, 
         \right)^{\frac{p+2r+1}{3p+1}}
         \, \mathrm{d}\tau \\[5pt]
& \leq \left\{\begin{array} {ll} 
      C( \, \alpha, \, \epsilon, \, p, \, r, \, \phi_0, \, \partial _x u_0 \, ) \, 
      \displaystyle{ \int ^t_0 
      (1+\tau )
      ^{\alpha - \frac{2p+r+1}{3p+1} -{\frac{p(p+2r+1)}{(p+1)(3p+1)}}} 
      \, \mathrm{d}\tau } \\[15pt] 
      \, \: \: \; \; \quad \quad \qquad \qquad \qquad \, 
      \left( \, 1 < p \le 
      \displaystyle{\frac{7}{12} + 
               \sqrt{\frac{73}{144} - \frac{(p+1)(3p-2)}{3}\, \epsilon } }
      \, \right),\\[15pt] 
      C( \, \alpha, \, \epsilon, \, p, \, r, \, \phi_0, \, \partial _x u_0 \, ) \, 
      \displaystyle{ \int ^t_0 
      (1+\tau )
      ^{\alpha - \frac{2p+r+1}{3p+1} -{\frac{3(p+2r+1)}{2(3p+1)(3p-2)}}
      + \frac{p+2r+1}{3p+1} \epsilon} 
      \, \mathrm{d}\tau } \\[15pt] 
      \, \: \: \; \; \quad \quad \qquad \qquad \qquad \qquad 
      \left( \,  \displaystyle{\frac{7}{12} + 
               \sqrt{\frac{73}{144} - \frac{(p+1)(3p-2)}{3}\, \epsilon } }
      < p \, \right)
      \end{array}
      \right.\, \\[5pt] 
& \leq \left\{\begin{array} {ll} 
      C( \, \alpha, \, \epsilon, \, p, \, r, \, \phi_0, \, \partial _x u_0 \, ) \, 
      (1+t )
      ^{\alpha - \frac{r}{p+1} } \\[15pt] 
      \, \: \: \; \; \quad \quad \qquad \qquad \qquad \, 
      \left( \, 1 < p \le 
      \displaystyle{\frac{7}{12} + 
               \sqrt{\frac{73}{144} - \frac{(p+1)(3p-2)}{3}\, \epsilon } }
      \, \right),\\[15pt] 
      C( \, \alpha, \, \epsilon, \, p, \, r, \, \phi_0, \, \partial _x u_0 \, ) \, 
      (1+t )
      ^{\alpha - \frac{6p(r-p)+7p+2r+3}{2(3p+1)(3p-2)} 
      + \frac{p+2r+1}{3p+1} \epsilon} \\[15pt] 
      \, \: \: \; \; \quad \quad \qquad \qquad \qquad \qquad 
      \left( \,  \displaystyle{\frac{7}{12} + 
               \sqrt{\frac{73}{144} - \frac{(p+1)(3p-2)}{3}\, \epsilon } }
      < p \, \right),
      \end{array}
      \right.\, 
\end{aligned}
\end{align}
\begin{align}
\begin{aligned}
&C_{p,r} \int ^t_0 (1+\tau )^\alpha 
    \left( \, 
    \int _{\partial _x u < 0 } 
    \left|\, \partial _x u \, \right|^{p+1} 
    \, \mathrm{d}x \, 
    \right)^{\frac{p+2r+2}{3p}}
    \, \mathrm{d}\tau  \\
& \leq C_{p,r} \int ^t_0 (1+\tau )^\alpha 
       \left( \, 
       \int _{-\infty}^{\infty}
       \left|\, \partial _x \phi \, \right|^{p+1} 
       \, \mathrm{d}x \, 
       \right)
       \| \, \partial _x u (\tau) \, \| _{L^{p+1}}
       ^{\frac{2(p+1)(r-p+1)}{3p}} 
       \, \mathrm{d}\tau \\[5pt]
& \leq \left\{\begin{array} {ll} 
      C( \, \alpha, \, \epsilon, \, p, \, r, \, \phi_0, \, \partial _x u_0 \, ) \, 
      \displaystyle{ \int ^t_0 
      (1+\tau )
      ^{\alpha - \frac{2(r-p+1)}{3(p+1)} } 
      \| \, \partial _x \phi (\tau) \, \| _{L^{p+1}}^{p+1}
      \, \mathrm{d}\tau } \\[15pt] 
      \, \: \: \; \; \; \; \; \quad \quad \qquad \qquad \qquad \, 
      \left( \, 1 < p \le 
      \displaystyle{\frac{7}{12} + 
               \sqrt{\frac{73}{144} - \frac{(p+1)(3p-2)}{3}\, \epsilon } }
      \, \right),\\[15pt] 
      C( \, \alpha, \, \epsilon, \, p, \, r, \, \phi_0, \, \partial _x u_0 \, ) \, 
      \displaystyle{ \int ^t_0 
      (1+\tau )
      ^{\alpha - \frac{r-p+1}{p(3p-2)} + \frac{2(r-p+1)}{3p} \epsilon} 
      \| \, \partial _x \phi (\tau) \, \| _{L^{p+1}}^{p+1}
      \, \mathrm{d}\tau } \\[15pt] 
      \, \: \: \; \; \; \; \; \quad \quad \qquad \qquad \qquad \qquad 
      \left( \,  \displaystyle{\frac{7}{12} + 
               \sqrt{\frac{73}{144} - \frac{(p+1)(3p-2)}{3}\, \epsilon } }
      < p \, \right)
      \end{array}
      \right.\, 
\end{aligned}
\end{align}
for any $0<\epsilon \ll 1$. \\
\noindent
By using Lemma 5.9 with 
\begin{align*}
\begin{aligned}
\alpha \mapsto \left\{\begin{array} {ll} 
               \alpha-\displaystyle{\frac{2(r-p+1)}{3(p+1)}} 
               \, \, \: \; \quad \quad \qquad \qquad 
               \left( \, 1 < p \le 
               \displaystyle{\frac{7}{12} + 
               \sqrt{\frac{73}{144} - \frac{(p+1)(3p-2)}{3}\, \epsilon } }
               \, \right),\\[15pt] 
               \alpha
               -\displaystyle{\left( \, 
                \frac{r-p+1}{p(3p-2)} - \frac{2(r-p+1)}{3p} \epsilon 
                \, \right)}
               \quad 
               \left( \,  \displaystyle{\frac{7}{12} + 
               \sqrt{\frac{73}{144} - \frac{(p+1)(3p-2)}{3}\, \epsilon } }
               < p \, \right)
               \end{array}
               \right.\, 
\end{aligned}
\end{align*}
and $q=2$, 
we get 
\begin{align}
\begin{aligned}
&C_{p,r} \int ^t_0 (1+\tau )^\alpha 
    \left( \, 
    \int _{\partial _x u < 0 } 
    \left|\, \partial _x u \, \right|^{p+1} 
    \, \mathrm{d}x \, 
    \right)^{\frac{p+2r+2}{3p}}
    \, \mathrm{d}\tau  \\[5pt]
& \leq \left\{\begin{array} {ll} 
      C( \, \alpha, \, \epsilon, \, p, \, r, \, \phi_0, \, \partial _x u_0 \, ) \, 
      (1+t )
      ^{\alpha - \frac{4p(r-p)+7p+3}{6p} } \\[15pt] 
      \, \: \: \; \; \quad \qquad \qquad \qquad \qquad \, 
      \left( \, 1 < p \le 
      \displaystyle{\frac{7}{12} + 
        \sqrt{\frac{73}{144} - \frac{(p+1)(3p-2)}{3}\, \epsilon } }
      \, \right),\\[15pt] 
      C( \, \alpha, \, \epsilon, \, p, \, r, \, \phi_0, \, \partial _x u_0 \, ) \, 
      (1+t )
      ^{\alpha - \frac{p+2r}{2p(3p-2)} + \frac{2(r-p+1)}{3p} \epsilon} \\[15pt] 
      \, \: \: \; \; \quad \qquad \qquad \qquad \qquad \qquad 
      \left( \,  \displaystyle{\frac{7}{12} + 
        \sqrt{\frac{73}{144} - \frac{(p+1)(3p-2)}{3}\, \epsilon } }
      < p \, \right)
      \end{array}
      \right.\, 
\end{aligned}
\end{align}
for any $0<\epsilon \ll 1$. \\
\noindent
Substituting (6.3) and (6.5) into (6.2), 
we have 
\begin{align}
\begin{aligned}
 &(1+t)^\alpha  
  \| \, \partial _x u (t)\, \| _{L^{r+1}}^{r+1} \\
 & \quad 
  + \int ^t_0 (1+\tau )^\alpha  \int _{-\infty}^{\infty} 
    \bigl| \, \partial _x u \, \bigr|^{p+r-2} 
    \left( \, \partial _x^2 u \, \right)^2 
    \, \mathrm{d}x \mathrm{d}\tau \\
 & \quad 
  + \int ^t_0 (1+\tau )^\alpha  
    \int _{\partial _x u \geq 0} 
    f''(u) \left|\, \partial _x u \, \right|^{r+2} 
    \, \mathrm{d}x \mathrm{d}\tau  \\
& \leq C_{\alpha,p,r} 
       \| \, \partial _x u_{0} \, \| _{L^{p+1}} \\[5pt]
& \quad + 
      \left\{\begin{array} {ll} 
      C( \, \alpha, \, \epsilon, \, p, \, r, \, \phi_0, \, \partial _x u_0 \, ) 
      (1+t)^{\alpha}\\[5pt]
      \times 
      \left( \, 
       (1+t)
       ^{- \frac{r}{p+1}} 
       + (1+t)
       ^{- \frac{4p(r-p)+7p+3}{6p(p+1)} } 
       \,  \right) \\[15pt] 
      \, \: \: \; \; \quad \qquad \qquad \qquad \, 
      \left( \, 1 < p \le 
      \displaystyle{\frac{7}{12} + 
      \sqrt{\frac{73}{144} - \frac{(p+1)(3p-2)}{3}\, \epsilon } }
      \, \right),\\[15pt] 
      C( \, \alpha, \, \epsilon, \, p, \, r, \, \phi_0, \, \partial _x u_0 \, ) 
      (1+t)^{\alpha}\\[5pt]
      \times 
      \left( \, 
       (1+t)
       ^{- \frac{6p(r-p)+7p+2r+3}{2(3p+1)(3p-2)} 
       + \frac{p+2r+1}{3p+1} \epsilon} 
       + (1+t)
       ^{- \frac{p+2r}{2p(3p-2)} 
       + \frac{2(r-p+1)}{3p} \epsilon} 
       \,  \right) \\[15pt] 
      \, \: \: \; \; \quad \qquad \qquad \qquad \qquad 
      \left( \,  \displaystyle{\frac{7}{12} + 
      \sqrt{\frac{73}{144} - \frac{(p+1)(3p-2)}{3}\, \epsilon } }
      < p \, \right)
      \end{array}
      \right.\, \\[5pt] 
\end{aligned}
\end{align}
for any $0<\epsilon \ll 1$.\\
\noindent
Here, we note the following: 
if $1 < p 
    \le \frac{7}{12} + 
        \sqrt{\frac{73}{144} - \frac{(p+1)(3p-2)}{3}\, \epsilon }
    \left( \, < \frac{7+\sqrt{73}}{12} \, \right)$, 
then 
$$
p<\frac{-4p^2+7p+3}{2p}. 
$$
Therefore, it follows that 
$$
(1+t)^{- \frac{r}{p+1} } 
\leq (1+t)^{- \frac{4p(r-p)+7p+3}{6p(p+1)}} \quad 
\left( \, r>\frac{-4p^2+7p+3}{2p}>p \, \right)
$$
and 
$$
(1+t)^{- \frac{4p(r-p)+7p+3}{6p(p+1)}}
\leq (1+t)^{- \frac{r}{p+1} }  \quad 
\left( \, p<r\leq \frac{-4p^2+7p+3}{2p} \, \right).
$$
If $\epsilon 
< \frac{3(r-p)(p-1)(3p+1)}{(3p-2)|\, 9p^2-p-2r-2\, |}$, 
then, for any $p>1$ with $r>p$, 
$$
(1+t)^{- \frac{6p(r-p)+7p+2r+3}{2(3p+1)(3p-2)} 
       + \frac{p+2r+1}{3p+1} \epsilon} 
\leq (1+t)^{- \frac{p+2r}{2p(3p-2)} 
       + \frac{2(r-p+1)}{3p} \epsilon}. 
$$

\medskip

\noindent
Thus, we do complete the proof of Proposition 6.1.

%
%








\bibliographystyle{model6-num-names}
\bibliography{<your-bib-database>}







\end{document}